\documentclass[11pt, reqno]{amsart}
\usepackage{amsmath, amsthm, amscd, amsfonts, amssymb, graphicx, color}
\usepackage[bookmarksnumbered, colorlinks, plainpages]{hyperref}

\newtheorem{theorem}{Theorem}[section]
\newtheorem{lemma}[theorem]{Lemma}
\newtheorem{proposition}[theorem]{Proposition}

\theoremstyle{definition}

\theoremstyle{remark}

\numberwithin{equation}{section}

\input{mathrsfs.sty}

\begin{document}
\setcounter{page}{1}

\title{Essential Normality of automorphic composition
operators}

\author{Liangying Jiang, Caiheng Ouyang,  Ruhan Zhao}

\address{Department of Applied Mathematics,  Shanghai Finance University, Shanghai 201209, P. R. China}

\email{\textcolor[rgb]{0.00,0.00,0.84}{liangying1231@163.com,
jiangly@shfc.edu.cn}}

\address{Wuhan Institute of Physics and Mathematics, Chinese Academy of Sciences,
Wuhan 430071, P. R. China}

\email{\textcolor[rgb]{0.00,0.00,0.84}{ ouyang@wipm.ac.cn}}

\address{ Department of  Mathematics, The College at Brockport, State University of New York, Brockport, New York 14420, USA}
\email{\textcolor[rgb]{0.00,0.00,0.84}{ rzhao@brockport.edu}}



\date{March 5th, 2014
\newline \indent Supported by the National
Natural Science Foundation of China (No.11101279 and No.11271359)}

\begin{abstract} We first characterize those  composition operators that are essentially normal on the weighted Bergman space  $A^2_s(D)$ for any real $s>-1$,
where induced symbols are automorphisms of the unit  disk $D$. Using the  same technique, we investigate the automorphic composition operators on  the Hardy space   $H^2(B_N)$  and   the weighted Bergman spaces $A^2_s(B_N)$ ($s>-1$). Furthermore, we give some  composition operators induced by
 linear fractional self-maps of the unit ball $B_N$ that are not essentially normal.
\end{abstract} \maketitle

\section{Introduction and preliminaries}

 Let $B_N=\{z\in \mathbb{C}^N: |z|<1 \}$ denote the unit ball
of  $\mathbb{C}^N$ and  $\partial B_N$ denote  the boundary of $B_N$. The Hardy space $H^2(B_N)$  consists  of holomorphic
functions $f$ in $B_N$  such that   $$||f||^2\equiv
\sup\limits_{0<r<1}\int_{\partial B_N}|f(r\zeta)|^2
d\sigma(\zeta)<\infty,$$
where $d\sigma$ denotes the normalized surface measure on $\partial B_N$. Let  $dv$ denote the normalized volume measure on $B_N$.
For $s>-1$, the weighted Lebesgue measure $$dv_s(z)=c_s(1-|z|^2)^sdv(z):=\frac{\Gamma(N+s+1)}{N!\Gamma(s+1)}(1-|z|^2)^sdv(z).$$
 The weighted Bergman space $A^2_s(B_N)$ consists of holomorphic functions $f$ in $B_N$ satisfying
$$ ||f||^2_s\equiv
\int_{B_N}|f(z)|^2 dv_s(z)<\infty.$$  In this article, we let $\mathcal{H}$ denote the Hardy space $H^2(B_N)$ or the weighted Bergman space $A^2_s(B_N)$ ($s>-1$).

Let $\varphi$ be a holomorphic self-map of $B_N$, the
 composition operator $C_\varphi$ on the space $\mathcal{H}$ is defined by   $$C_\varphi(f)=f\circ\varphi$$ for $f\in\mathcal{H}$.

An operator $T$ on the space $H$ is essentially normal if its self-commutator $[T^\ast, T]=T^\ast T-T T^\ast$ is compact. Equivalently, an operator is essentially normal if its image in the Calkin algebra $B(H)/B_0(H)$ is normal.  A surprising  result is that essentially normal operators can be  characterized up to unitary equivalence modulo the compact operators.

In this article,  we are interested in which  composition operators are essentially normal on some classical Hilbert spaces. This question is difficult to  answer even on the Hardy space $H^2(D)$, unless $\varphi$ is an automorphism or $\varphi$ is a linear fractional self-map of the unit disk $D$.

When  $\varphi$ is an automorphism of $D$, the operator $C_\varphi$ is essentially normal on $H^2(D)$ if and only if $\varphi$ is a rotation, i.e. $C_\varphi$ is normal (see \cite{BLNS}). In \cite{MW04}, this result was extended to the weighted Bergman space $A^2_s(D)$  for any positive integer $s$. Soon after, MacCluer and Pons \cite{MP06} obtained the same result for $C_\varphi$ acting on $H^2(B_N)$ and the weighted Bergman space $A^2_s(B_N)$ for real $s>-1$, where $\varphi$ is an automorphism of $B_N$. Moreover, on  the Hardy space $H^2(B_N)$ and the Bergman space $A^2(B_N)$, a very simple proof of this result can be found in \cite{BM07}.

In Section 2 of this paper, we get
\\ \\
{\bf Theorem 2.4.} {\it If $\varphi$ is a non-rotation automorphism of $D$, then $C_\varphi$  is not essentially normal
on $A^2_s(D)$ for any real $s>-1$.}
\\ \\ This means that $C_\varphi$ is essentially normal on $A^2_s(D)$ for any real $s>-1$ if and only if $\varphi$ is a rotation.
For the case of the unit ball, we have the following result in Section 3.
\\ \\
{\bf Theorem 3.7.} {\it Let $\varphi$ be  an automorphism  of $B_N$ with  $a=\varphi^{-1}(0)\ne 0$ and
 $$e_m(z)=\sqrt{\frac{\Gamma(t+m)}{\Gamma(t)!m!}}\biggl<z, \frac{a}{|a|}\biggr>^m.$$
Then on the space $\mathcal{H}$,
$$\lim\limits_{m\to\infty}(||C^*_\varphi(e_m)||_{\mathcal{H}}^2-||C_\varphi(e_m)||_{\mathcal{H}}^2)>0.$$
Here $t=N$ when $\mathcal{H}=H^2(B_N)$ and $t=N+s+1$ when $\mathcal{H}=A^2_s(B_N)$ ($s>-1$).}
\\ \\
Therefore, as an immediate result of this theorem, we see that  $C_\varphi$  is essentially normal on the Hardy space $H^2(B_N)$ or the  weighted  Bergman space $A_s^2(B_N)$ ($s>-1$) if and only if $\varphi$ is unitary. Theorem 3.7 will also be used for the discussion about linear fractional composition operators in Section 4.

If $\varphi$ is a linear fractional self-map of $D$, using the adjoint formula $C^*_\varphi=T_gC_\sigma T^*_h$ (Here, $g, \sigma, h$ will be introduced in Section 2), we have known that $C_\varphi$ is essentially normal on $H^2(D)$ or  $A^2_s(D)$  ($s>-1$) if and only if $\varphi$ is a parabolic non-automorphism (see \cite{BLNS} and \cite{MW04}).

For linear fractional self-maps  of $B_N$, the situation is more complicated. Especially, from the discussion of  spectral  structures of linear fractional composition operators on $H^2(B_N)$ (see \cite{B10}, \cite{BC12}, \cite{JC14}), we have found that  linear fractional maps of $B_N$,  conjugated by automorphisms,  must be classified into nine different cases.  Until now, we only know a little  about which linear fractional composition operators are essentially normal on the space $\mathcal{H}$ (see  \cite{JO09} and \cite{MW05}). Even when $\varphi$ is  positive parabolic, i.e. $\varphi$ is parabolic with $\varphi\circ\sigma=\sigma\circ\varphi$, where $\sigma$ is the adjoint map of $\varphi$, we do not know whether $C_\varphi$ is  essentially normal  on $\mathcal{H}$.

In Section 4, we are interested in the essential normality of $C_\varphi$ on $H^2(B_N)$ and  $A_s^2(B_N)$ ($s>-1$) when $\varphi$  is a linear fractional self-map of $B_N$. The main results are as follows.
\\ \\
{\bf Theorem 4.1.} {\it If  $\varphi$ is   a linear fractional self-map   of $B_N$ and  the restriction of $\varphi$ to
$$B_k=\{(z_1,\ldots, z_N)\in B_N;\  z_i=0\ \mbox{for}\ i> k\}$$
 is a non-rotation automorphism of $B_k$,  then $C_\varphi$ is not essentially normal on $H^2(B_N)$ or $A^2_s(B_N)$ for any real $s>-1$.}
\\ \\
{\bf Theorem 4.4.} {\it Let  $\varphi$ be  a linear fractional self-map   of $B_N$ with only one interior fixed point $z_0$ on $\overline{B_N}$.
 If $p=\mbox{dim}\,L_U(\varphi, z_0)=0$ and  $||\varphi||_{\infty}=1$, then $C_\varphi$ is not essentially normal on  $H^2(B_N)$ or $A^2_s(B_N)$ ($s>-1$).  Where $L_U(\varphi, z_0)$ is the unitary space of $\varphi$ at $z_0$ (see  Definition 2).}
\\ \par
One may be   surprised why we pay attention to the essential normality of $C_\varphi$ again, when $\varphi$ is an automorphism of $D$ or $B_N$. We consider this problem based on two reasons: First, we want to exhibit how we combine perfectly the idea of MacCluer and Weir \cite{MW04}  with a tool provided by MacCluer and Pons \cite{MP06}. Moreover, on the unit ball, an astonishing result is that the sequence $\{<z, \zeta>^n\}$ for $\zeta\in\partial B_N$ behaves more similarly to the basis $\{z^n\}$ of $D$  (Also, one may compare the condition for composition operators to be  compact on the Bloch space $\mathcal{B}(B_N)$ given by Dai \cite{Dai12}). Second,  using  Theorem 3.7  obtained in Section 3,  we show that some linear fractional composition operators are not essentially normal on the space  $\mathcal{H}$.


\section { Essential normality of  composition operators in the unit disk}

In this section, we consider the essential normality of composition operators   induced by  automorphisms of $D$ on the weighted Bergman space $A^2_s(D)$ for any real $s>-1$.

Recall that  a  linear fractional map of $B_N$ is of the form
$$\varphi(z)=\frac{Az+B}{<z, C>+d}$$  with $A\in \mathbb{C}^{N\times N}$, $B, C\in \mathbb{C}^{N\times 1}$, $d\in
 \mathbb{C}$, where  $<\cdot,
\cdot>$ denotes the Euclidean inner product on $\mathbb{C}^N$. If  $\varphi$ is a  linear
fractional self-map  of $B_N$, we have the following
adjoint formula for $C_\varphi$ on the space $\mathcal{H}$ ( see
\cite{CM00} or \cite{MW05}),
$$C_{\varphi}^{\ast}=T_g C_{\sigma} T_h^{\ast},$$
where $$\sigma(z)=\frac{A^\ast z-C}{<z,  -B>+ \bar{d}}$$ is the
adjoint map of $\varphi$,  $T_g$ and $T_h$ are  analytic Toeplitz
operators respectively with  symbols $g(z)=(<z,  -B>+ \bar{d})^{-t}$  and
$h(z)=(<z, C>+ d)^t$.  Here, $t=N$ when $\mathcal{H}=H^2(B_N)$  and $t=N+s+1$
when  $\mathcal{H}=A^2_s(B_N)$ ($s>-1$). On $H^2(D)$, this adjoint formula is often called Cowen's adjoint formula (see  \cite{Co88}). For the case of  $A^2_s(D)$ ($s>-1$) see  the reference \cite{Hu97}.

In order to determine which automorphic composition operators are essentially normal on the Bergman space $A^2(D)$, MacCluer and Weir \cite{MW04} applied the adjoint formula $C^*_\varphi=T_gC_\sigma T^*_h$  to calculate $\lim\limits_{n\to\infty}||C_\varphi^*(e_n)||_{A^2(D)}^2$ with $e_n(z)=\sqrt{n+1}z^n$.  In  \cite{MW04}, they also computed
$\lim\limits_{n\to\infty}||C_\varphi(e_n)||_{A^2(D)}^2$ and compared the two limits to deduce a necessary condition for  $C_\varphi$ to be essentially normal.
On the weighted Bergman space $A^2_s(D)$ for any positive integer $s$, they mentioned that, without proof, similar computations can give the following limits $$\lim\limits_{n\to\infty}||C_\varphi(e_n)||_s^2 \qquad \mbox{and} \qquad \lim\limits_{n\to\infty}||C^*_\varphi(e_n)||_s^2,$$ where $$e_n(z)=\sqrt{\frac{\Gamma(s+2+n)}{\Gamma(s+2)\Gamma(n+1)}}z^n.$$

But it is surprising that we can get the same  limit $\lim\limits_{n\to\infty}||C^*_\varphi(e_n)||_s^2$ on $A^2_s(D)$ for any real $s>-1$  using the following formula (see \cite{MP06}) $$C_{\varphi}C_{\varphi}^*=T_f+K, \eqno(2.1)$$
where  $\varphi$ is an automorphism  of $B_N$,  $T_f$ is the Toeplitz operator with symbol $f$ and $K$ is a compact operator on $\mathcal{H}$. Next, we will give a detail computation for  $\lim\limits_{n\to\infty}||C^*_\varphi(e_n)||_s^2$.

First, we need describe exactly the symbol $f$ in the formula (2.1) when $\varphi=\varphi_a$ is an involution automorphism of $B_N$ which interchanges $a$ and $0$.

\begin{proposition}\label{prop 2.1}
Let $\varphi_a(z)$ be an involution automorphism of $B_N$, interchanging $a$ and $0$. Then
$$C_{\varphi_a}C_{\varphi_a}^*=T_f+K,$$
where $T_f$ is the Toeplitz operator with symbol  $$f(z)=\biggl(\frac{|1-<z,a>|^2}{1-|a|^2}\biggr)^t$$
and $K$ is a compact operator on $\mathcal{H}$. Here, $t=N$ when $\mathcal{H}=H^2(B_N)$  and $t=N+s+1$ when $\mathcal{H}=A^2_s(B_N)$ ($s>-1$).
\end{proposition}
\proof This result can be easily obtained by Proposition 1 in \cite{MP06}. For completeness, we give a simple proof.

Recall that $$\varphi_a(z)=\frac{a-P_a(z)-s_aQ_a(z)}{1-<z,a>}, \quad z\in B_N$$ for any point $a\in B_N-\{0\}$ and $\varphi_0(z)=-z$,
where $s_a=\sqrt{1-|a|^2}$, $$P_a(z)=\frac{<z,a>}{|a|^2}a \qquad \mbox{and}\qquad Q_a(z)=z-\frac{<z,a>}{|a|^2}a.$$ Thus, $C_{\varphi_a}^*=T_gC_\sigma T_h^*$
 with $g(z)=(1-<z,a>)^{-t}$  and $h(z)=(1-<z,a>)^t$, where $t=N$ when $\mathcal{H}=H^2(B_N)$  and $t=N+s+1$ when $\mathcal{H}=A^2_s(B_N)$. Moreover, we have $\sigma=\varphi_a^{-1}=\varphi_a$ by Lemma 6.3 of \cite{BB02}. Using the semi-multiplicative property for Toeplitz operator mod $\mathcal{K}$ (see \cite{MP06} for details), we see that
$$C_{\varphi_a}C_{\varphi_a}^*=C_{\varphi_a}T_gC_\sigma T_h^*=T_{g\circ \varphi_a}C_{\varphi_a}C_\sigma T_h^*=T_{g\circ \varphi_a} T_h^*\equiv T_{\overline{h}\, (g\circ \varphi_a)}\ (\mbox{mod}\, \mathcal{K}),$$ where $\mathcal{K}$ denotes the ideal  of compact operators on $\mathcal{H}$.
Now, applying the equality $$1-<\varphi_a(z),a>=\frac{1-|a|^2}{1-<z,a>},$$ we get  $$f(z)=\overline{h}(z) (g\circ \varphi_a(z))=\biggl(\frac{\overline{1-<z,a>}}{1-<\varphi_a(z),a>}\biggr)^t=\biggl(\frac{|1-<z,a>|^2}{1-|a|^2}\biggr)^t.$$
This gives the desired conclusion.
\ \ $\Box$
\\ \par
Now, we will show that how the formula in Proposition 2.1 works on the calculation of  $\lim\limits_{n\to\infty}||C_\varphi^*(e_n)||_{A^2(D)}^2$. After that, we will use similar method to compute  the corresponding limit  on $A^2_s(D)$ for any real $s>-1$.

\begin{proposition}\label{prop 2.2}
If $\varphi$ is an automorphism of $D$  and $e_n(z)=\sqrt{n+1}z^n$, then on $A^2(D)$,
$$\lim\limits_{n\to\infty}||C_\varphi^*(e_n)||_{A^2(D)}^2=\frac{1+4|p|^2+|p|^4}{(1-|p|^2)^2},$$ where $p=\varphi^{-1}(0)$.
\end{proposition}
\proof This is Proposition 2 in \cite{MW04}, we will use Proposition 2.1 to give another proof.

Any automorphism $\varphi$ of $D$ with  $\varphi(p)=0$ can be written as $\varphi=U\circ \varphi_p$, where $U(z)=\lambda z$ with $|\lambda|=1$
and $\varphi_p(z)=\frac{p-z}{1-\overline{p}z}$. Thus, $C_\varphi=C_{U\circ \varphi_p}=C_{\varphi_p}C_U$. Applying Proposition 2.1, we see that
$$C_\varphi C_\varphi^*=C_{\varphi_p}C_U(C_{\varphi_p}C_U)^*=C_{\varphi_p}C_U C_U^*C_{\varphi_p}^*=C_{\varphi_p}C_{\varphi_p}^*=T_f+K$$
with $$f(z)=\biggl(\frac{|1-\overline{p}z|^2}{1-|p|^2}\biggr)^2$$ and $K$ is a compact operator on $A^2(D)$. This implies that $$\lim\limits_{n\to\infty}||C_\varphi^*(e_n)||_{A^2(D)}^2=\lim\limits_{n\to\infty}<C_\varphi C_\varphi^*(e_n), e_n>=\lim\limits_{n\to\infty}<T_f(e_n), e_n>.$$

Write $$f(z)=\frac{1}{(1-|p|^2)^2}(1-\overline{p}z)^2(1-p\overline{z})^2,$$ then $$T_f(e_n)=\frac{\sqrt{n+1}}{(1-|p|^2)^2}P[(1-\overline{p}z)^2(1-p\overline{z})^2 z^n],$$ where $P$ denotes the orthogonal projection of $L^2(D)$
onto $A^2(D)$. Note that $$P(\overline{z} z^n)=T^*_z z^n=\frac{n}{n+1}z^{n-1} \quad \mbox{for}\ n\ge 1$$
 and $$P(\overline{z}^2 z^n)=T^*_{z^2} z^n=\frac{n-1}{n+1}z^{n-2} \quad \mbox{for}\ n\ge 2.$$ We compute that{\setlength\arraycolsep{2pt}
\begin{eqnarray*}&& P[(1-\overline{p}z)^2(1-p\overline{z})^2 z^n]
 \\ &=& P(z^n-2p\overline{z} z^n+p^2\overline{z}^2 z^n-2\overline{p} z^{n+1}+4|p|^2\overline{z} z^{n+1}-2|p|^2p\overline{z}^2 z^{n+1}
 \\ && +\overline{p}^2 z^{n+2}-2|p|^2\overline{p}\cdot\overline{z}  z^{n+2}+|p|^4\overline{z}^2  z^{n+2})
\\ &=& z^n-2\frac{n}{n+1}pz^{n-1}+\frac{n-1}{n+1}p^2 z^{n-2}-2\overline{p}z^{n+1}+4\frac{n+1}{n+2}|p|^2 z^n
\\ && -2\frac{n}{n+2}|p|^2p z^{n-1}+\overline{p}^2  z^{n+2}- 2\frac{n+2}{n+3}|p|^2\overline{p} z^{n+1}+\frac{n+1}{n+3} |p|^4 z^n.
\end{eqnarray*}}It follows that{\setlength\arraycolsep{2pt}\begin{eqnarray*} <T_f(e_n), e_n>&=& \frac{n+1}{(1-|p|^2)^2}<P[(1-\overline{p}z)^2(1-p\overline{z})^2 z^n], z^n>
 \\ &=&\frac{n+1}{(1-|p|^2)^2}\biggr(||z^n||^2+4\frac{n+1}{n+2}|p|^2||z^n||^2+\frac{n+1}{n+3} |p|^4||z^n||^2\biggl).
\end{eqnarray*}}Since $||z^n||^2:=||z^n||^2_{A^2(D)}=\frac{1}{n+1}$,  taking the limit, we obtain
{\setlength\arraycolsep{2pt}\begin{eqnarray*}\lim\limits_{n\to\infty}||C_\varphi^*(e_n)||_{A^2(D)}^2 &=& \lim\limits_{n\to\infty}<T_f(e_n), e_n>
 \\ &=& \lim\limits_{n\to\infty}\frac{1}{(1-|p|^2)^2}\biggr(1+4\frac{n+1}{n+2}|p|^2+\frac{n+1}{n+3} |p|^4\biggl)
  \\ &=&\frac{1+4|p|^2+|p|^4}{(1-|p|^2)^2}.
\end{eqnarray*}}
 $\Box$
\\ \par
Using this technique, we easily get  the following result on the weighted Bergman space $A_s^2(D)$ for real $s>-1$. First, we need some notations. For any real $c$, we denote \begin{displaymath}
\left(\begin{array}{c}
c
\\ 0
\end{array}\right)=1 \quad \mbox{and} \quad
\left(\begin{array}{c}
c
\\ k
\end{array}\right)=\frac{c(c-1)\cdots(c-k+1)}{k!}, \ k\ge 1.
\end{displaymath}
Let $(c)_k$ denote the shifted factorial defined by $$(c)_k=c(c+1)\cdots(c+k+1)=\frac{\Gamma(c+k)}{\Gamma(c)}\quad
\mbox{for}\ k>0, \ (c)_0=1,$$ where $c$ is any real or complex number. Write
$$F((a_1)_k,\cdots,(a_p)_k;(b_1)_k,\cdots,(b_q)_k;x)=\sum\limits^{\infty}_{k=0}\frac{(a_1)_k\cdots(a_p)_k}{(b_1)_k\cdots(b_q)_k}\frac{x^k}{k!},$$
which is a hypergeometric function (see \cite{AAR}).

\begin{proposition}\label{prop 2.3}
If $\varphi$ is an automorphism of $D$  and $$e_n(z)=\sqrt{\frac{\Gamma(s+2+n)}{\Gamma(s+2)\Gamma(n+1)}}z^n,$$ then on $A_s^2(D)$ for any real $s>-1$,
$$\lim\limits_{n\to\infty}||C_\varphi^*(e_n)||_s^2=\frac{1}{(1-|p|^2)^{s+2}}\sum\limits^{\infty}_{k=0}\left(\begin{array}{c}
s+2
\\ k
\end{array}\right)^2|p|^{2k},$$ where $p=\varphi^{-1}(0)$.
\end{proposition}
\proof First on $A_s^2(D)$ ($s>-1$), by Proposition 2.1,
{\setlength\arraycolsep{2pt}\begin{eqnarray*}f(z)&=&\biggl(\frac{|1-\overline{p}z|^2}{1-|p|^2}\biggr)^{s+2}
\\ &=&\frac{1}{(1-|p|^2)^{s+2}}(1-\overline{p}z)^{s+2}(1-p\overline{z})^{s+2}
\\ &=& \frac{1}{(1-|p|^2)^{s+2}}\sum\limits^{\infty}_{k=0}\left(\begin{array}{c}
s+2
\\ k
\end{array}\right)(-\overline{p}z)^k\sum\limits^{\infty}_{j=0}\left(\begin{array}{c}
s+2
\\ j
\end{array}\right)(-p\overline{z})^j
\\ &=& \frac{1}{(1-|p|^2)^{s+2}}\sum\limits^{\infty}_{k,j=0}\left(\begin{array}{c}
s+2
\\ k
\end{array}\right)\left(\begin{array}{c}
s+2
\\ j
\end{array}\right)(-1)^{k+j}\overline{p}^kp^jz^k\overline{z}^j.
\end{eqnarray*}}On the other hand,
$$P(\overline{z}^mz^l)=T^*_{z^m}(z^l)=\frac{\Gamma(s+2+l-m)}{\Gamma(s+2+l)}\cdot\frac{l!}{(l-m)!}z^{l-m}$$ for non-negative integers $l\ge m$ and $0$ otherwise, where $P$ denotes the projection of $L^2_s(D)$ onto $A^2_s(D)$  for $s>-1$.  Combining this with the orthogonality of $z^{k_1}$ and $z^{k_2}$ when $k_1\ne k_2$, we calculate that{\setlength\arraycolsep{2pt}\begin{eqnarray*} && <T_f(e_n), e_n>=
\frac{1}{(1-|p|^2)^{s+2}}\frac{\Gamma(n+s+2)}{\Gamma(s+2)\Gamma(n+1)}\times
\\ &&\sum\limits^{\infty}_{k,j=0}\left(\begin{array}{c}
s+2
\\ k
\end{array}\right)\left(\begin{array}{c}
s+2
\\ j
\end{array}\right)(-1)^{k+j}\overline{p}^kp^j<P(\overline{z}^jz^{k+n}), z^n>
\\ &=& \frac{1}{(1-|p|^2)^{s+2}}\frac{\Gamma(n+s+2)}{\Gamma(s+2)\Gamma(n+1)} \sum\limits^{\infty}_{k,j=0}\left(\begin{array}{c}
s+2
\\ k
\end{array}\right)\left(\begin{array}{c}
s+2
\\ j
\end{array}\right)(-1)^{k+j}\overline{p}^kp^j \times
\\ &&\frac{\Gamma(s+2+n+k-j)}{\Gamma(s+2+n+k)}\cdot\frac{(n+k)!}{(n+k-j)!}<z^{n+k-j}, z^n>
\\ &=&\frac{1}{(1-|p|^2)^{s+2}}\frac{\Gamma(n+s+2)}{\Gamma(s+2)\Gamma(n+1)}\times
\\ &&\sum\limits^{\infty}_{k=0}\left(\begin{array}{c}
s+2
\\ k
\end{array}\right)^2|p|^{2k}\frac{\Gamma(s+2+n)}{\Gamma(s+2+k+n)}\cdot\frac{(k+n)!}{n!}||z^n||^2_s
\\ &=&\frac{1}{(1-|p|^2)^{s+2}}\sum\limits^{\infty}_{k=0}\left(\begin{array}{c}
s+2
\\ k
\end{array}\right)^2|p|^{2k}\frac{\Gamma(s+2+n)}{\Gamma(s+2+k+n)}\cdot\frac{(k+n)!}{n!},
\end{eqnarray*}}where in the last line we have used the norm $$||z^n||^2_s=\frac{\Gamma(s+2)\Gamma(n+1)}{\Gamma(n+s+2)}.$$

For any fixed non-negative integer $k$, using  Stirling's formula, we see that
$$\frac{\Gamma(s+2+n)}{\Gamma(s+2+k+n)}\cdot\frac{(k+n)!}{n!}\to 1$$ as $n\to\infty$. Hence, $|p|<1$ and the dominated convergence theorem give that{\setlength\arraycolsep{2pt}\begin{eqnarray*}
&&\lim\limits_{n\to\infty}<T_f(e_n), e_n>
\\ &&=\lim\limits_{n\to\infty}\frac{1}{(1-|p|^2)^{s+2}}\sum\limits^{\infty}_{k=0}\left(\begin{array}{c}
s+2
\\ k
\end{array}\right)^2|p|^{2k}\frac{\Gamma(s+2+n)}{\Gamma(s+2+k+n)}\cdot\frac{(k+n)!}{n!}
\\ &&=\frac{1}{(1-|p|^2)^{s+2}}\sum\limits^{\infty}_{k=0}\left(\begin{array}{c}
s+2
\\ k
\end{array}\right)^2|p|^{2k}.\end{eqnarray*}}Finally, using similar argument as in the proof of Proposition 2.2, we obtain
{\setlength\arraycolsep{2pt}\begin{eqnarray*} \lim\limits_{n\to\infty}||C_\varphi^*(e_n)||_s^2 & =& \lim\limits_{n\to\infty}<T_f(e_n), e_n>
\\ &=&\frac{1}{(1-|p|^2)^{s+2}}\sum\limits^{\infty}_{k=0}\left(\begin{array}{c}
s+2
\\ k
\end{array}\right)^2|p|^{2k}.\end{eqnarray*}}
 $\Box$
\\ \par
In \cite{MW04}, MacCluer and Weir  gave the following result  on the weighted Bergman space $A^2_s(D)$  only for any positive integer $s$. Now, using Proposition 2.3 and the idea of MacCluer and Weir (see Theorem 5 of \cite{MW04}), we can prove this result on $A^2_s(D)$ for any real $s>-1$.

\begin{theorem}\label{the 2.4} If $\varphi$ is a non-rotation automorphism of $D$, then $C_\varphi$  is not essentially normal
on $A^2_s(D)$ for any real $s>-1$.
\end{theorem}
\proof  Let $\{e_n\}$ be the normalized basis of $A^2_s(D)$ ($s>-1$) defined in Proposition 2.3.  Since  $\varphi$ is a non-rotation automorphism of $D$, we have
 $\varphi(z)=\lambda\frac{p-z}{1-\overline{p}z}$ with $|\lambda|=1$ and $p=\varphi^{-1}(0)\ne 0$. First, using the change of variables formula and the   orthogonality of $z^{k_1}$ and $z^{k_2}$ when $k_1\ne k_2$, we compute that  {\setlength\arraycolsep{2pt}
\begin{eqnarray*}&& ||C_\varphi(e_n)||^2_s = \frac{\Gamma(s+2+n)}{\Gamma(s+2)\Gamma(n+1)}||\varphi^n||^2_s
 \\ &=& \frac{\Gamma(s+2+n)}{\Gamma(s+2)\Gamma(n+1)}\int_D\biggl|\lambda\frac{p-z}{1-\overline{p}z}\biggr|^{2n} dv_s(z)
\\ &=&\frac{\Gamma(s+2+n)}{\Gamma(s+2)\Gamma(n+1)}\int_D |z|^{2n}\frac{(1-|p|^2)^{s+2}}{|1-\overline{p}z|^{2(s+2)}} dv_s(z)
\\ &=& (1-|p|^2)^{s+2}\frac{\Gamma(s+2+n)}{\Gamma(s+2)\Gamma(n+1)}\int_D |z|^{2n}\frac{1}{(1-\overline{p}z)^{s+2}}\frac{1}{(1-p\overline{z})^{s+2}} dv_s(z)
\\ &=&
(1-|p|^2)^{s+2}\frac{\Gamma(s+2+n)}{\Gamma(s+2)\Gamma(n+1)}\times
\\ && \sum\limits^{\infty}_{k,j=0}\frac{\Gamma(s+2+k)}{k!\Gamma(s+2)}\frac{\Gamma(s+2+j)}{j!\Gamma(s+2)}\overline{p}^k p^j\int_D |z|^{2n}z^k\overline{z}^jdv_s(z)
\\ & = &(1-|p|^2)^{s+2}\frac{\Gamma(s+2+n)}{\Gamma(s+2)\Gamma(n+1)}\sum\limits^{\infty}_{k=0}\biggl(\frac{\Gamma(s+2+k)}{k!\Gamma(s+2)}\biggr)^2|p|^{2k}\int_D |z|^{2(n+k)}dv_s(z)
\\ &=& (1-|p|^2)^{s+2}\sum\limits^{\infty}_{k=0}\biggl(\frac{\Gamma(s+2+k)}{k!\Gamma(s+2)}\biggr)^2
 |p|^{2k}\frac{\Gamma(s+2+n)}{\Gamma(n+1)}\cdot\frac{\Gamma(n+k+1)}{\Gamma(s+2+n+k)}.
\end{eqnarray*}}Since for any fixed non-negative integer $k$,
$$\frac{\Gamma(s+2+n)}{\Gamma(n+1)}\frac{\Gamma(n+k+1)}{\Gamma(s+2+n+k)}\to 1$$ as $n\to\infty$, as in the proof of Proposition 2.3, applying the dominated convergence theorem, we see that{\setlength\arraycolsep{2pt}
\begin{eqnarray*} && \lim\limits_{n\to\infty} ||C_\varphi(e_n)||^2_s
 \\ &=& \lim\limits_{n\to\infty}(1-|p|^2)^{s+2}\sum\limits^{\infty}_{k=0}\biggl(\frac{\Gamma(s+2+k)}{k!\Gamma(s+2)}\biggr)^2
 |p|^{2k}\frac{\Gamma(s+2+n)}{\Gamma(n+1)}\cdot\frac{\Gamma(n+k+1)}{\Gamma(s+2+n+k)}
 \\ &=& (1-|p|^2)^{s+2}\sum\limits^{\infty}_{k=0}\biggl(\frac{\Gamma(s+2+k)}{k!\Gamma(s+2)}\biggr)^2
 |p|^{2k}.
\end{eqnarray*}}

Next, we follow the same idea as in the proof of Theorem 5 in \cite{MW04} to deal with the above series. Using Euler's formula (see Theorem 2.25 of \cite{AAR}),
{\setlength\arraycolsep{2pt}
\begin{eqnarray*}&& \lim\limits_{n\to\infty} ||C_\varphi(e_n)||^2_s = (1-|p|^2)^{s+2}\sum\limits^{\infty}_{k=0}\biggl(\frac{\Gamma(s+2+k)}{k!\Gamma(s+2)}\biggr)^2
 |p|^{2k}
 \\ &=& (1-|p|^2)^{s+2}\sum\limits^{\infty}_{k=0}\frac{(s+2)_k(s+2)_k}{(1)_k}
 \frac{|p|^{2k}}{k!}
 \\ &=& (1-|p|^2)^{s+2}  F(s+2,s+2;1;|p|^2)
\\ &=& (1-|p|^2)^{s+2}(1-|p|^2)^{1-(s+2)-(s+2)}  F(1-(s+2),1-(s+2);1;|p|^2)
\\ &=& \frac{1}{(1-|p|^2)^{s+2}}(1-|p|^2) F(-s-1,-s-1;1;|p|^2)
\\ &=& \frac{1}{(1-|p|^2)^{s+2}}(1-|p|^2)\sum\limits^{\infty}_{k=0}\frac{(-s-1)_k(-s-1)_k}{(1)_k}
 \frac{|p|^{2k}}{k!}
\\ &=& \frac{1}{(1-|p|^2)^{s+2}}(1-|p|^2)\sum\limits^{\infty}_{k=0} \biggl(\frac{(s+1)s\cdots(s+2-k)}{k!}\biggr)^2|p|^{2k}.
\end{eqnarray*}}Note that{\setlength\arraycolsep{2pt}
\begin{eqnarray*} && (1-|p|^2)\sum\limits^{\infty}_{k=0} \biggl(\frac{(s+1)s\cdots(s+2-k)}{k!}\biggr)^2|p|^{2k}
\\ &=&\sum\limits^{\infty}_{k=0} \biggl(\frac{(s+1)s\cdots(s+2-k)}{k!}\biggr)^2|p|^{2k}-\sum\limits^{\infty}_{k=0} \biggl(\frac{(s+1)s\cdots(s+2-k)}{k!}\biggr)^2|p|^{2(k+1)}
\\ &=&\sum\limits^{\infty}_{k=0} \biggl(\frac{(s+1)s\cdots(s+2-k)}{k!}\biggr)^2|p|^{2k}-\sum\limits^{\infty}_{k=1} \biggl(\frac{(s+1)s\cdots(s+2-(k-1))}{(k-1)!}\biggr)^2|p|^{2k}
\\ &=&1+(s^2+2s)|p|^2+\sum\limits^{\infty}_{k=2} \biggl(\frac{(s+1)s\cdots(s+2-(k-1))}{k!}\biggr)^2[(s+2-k)^2-k^2]|p|^{2k}
\\ &=&1+(s^2+2s)|p|^2+\sum\limits^{\infty}_{k=2} \biggl(\frac{(s+1)s\cdots(s+2-(k-1))}{k!}\biggr)^2(s+2)(s+2-2k)|p|^{2k}.
\end{eqnarray*}}Therefore, since $p\ne 0$, by Proposition 2.3,
{\setlength\arraycolsep{2pt}
\begin{eqnarray*} &&  \lim\limits_{n\to\infty} (||C_\varphi^*(e_n)||^2_s-||C_\varphi(e_n)||^2_s)
\\ &=& \frac{1}{(1-|p|^2)^{s+2}}\biggl\{\sum\limits^{\infty}_{k=0}\biggl(\frac{(s+2)(s+1)\cdots(s+2-(k-1))}{k!}\biggr)^2|p|^{2k}
\\ && - \biggl[1+(s^2+2s)|p|^2 +\sum\limits^{\infty}_{k=2} \biggl(\frac{(s+1)s\cdots(s+2-(k-1))}{k!}\biggr)^2(s+2)(s+2-2k)|p|^{2k}\biggr]\biggr\}
\\ &=& \frac{1}{(1-|p|^2)^{s+2}}\biggl[(2s+4)|p|^2
\\ && +\sum\limits^{\infty}_{k=2} \biggl(\frac{(s+1)s\cdots(s+2-(k-1))}{k!}\biggr)^22k(s+2)|p|^{2k}\biggr]>0.
\end{eqnarray*}}

It is clear that $\{e_n\}$ is a weakly convergent sequence of $D$. However, if  $p=\varphi^{-1}(0)\ne 0$, we have shown that{\setlength\arraycolsep{2pt}
\begin{eqnarray*} \lim\limits_{n\to\infty} ||[C_\varphi^*,C_\varphi](e_n)||_s &\ge &  \lim\limits_{n\to\infty}|<[C_\varphi^*,C_\varphi](e_n), e_n>|
\\ &=&  \lim\limits_{n\to\infty} (||C_\varphi^*(e_n)||^2_s-||C_\varphi(e_n)||^2_s)>0.
\end{eqnarray*}}As a consequence, $[C_\varphi^*,C_\varphi]$ is not compact and hence $C_\varphi$ is not essentially normal on $A^2_s(D)$ for any real $s>-1$.
\ \  $\Box$

\section { Essentially  normality of  composition operators in the unit ball}
First, we introduce some notations. If $z=(z_1,\cdots,z_N)$ and  $w=(w_1,\cdots,w_N)$ are points in $\mathbb{C}^N$, we write $$<z,w>=z_1\overline{w_1}+\cdots+z_N\overline{w_N}\qquad \mbox{and} \qquad |z|=<z,z>^{1/2}.$$ For an $N$-tuple $\alpha=(\alpha_1,\ldots,\alpha_N)$ of non-negative integers, which is also  called a multi-index, we write $$|\alpha|=|\alpha_1|+\cdots+|\alpha_N|, \qquad \alpha!=\alpha_1!\cdots \alpha_N!$$ and $z^\alpha=z_1^{\alpha_1}\cdots z_N^{\alpha_N}$. If $\alpha$ and $\beta$ are two multi-indexes, we say $\beta\le \alpha$ provided $\beta_j\le \alpha_j$ for $1\le j\le N$. In this case, $\alpha-\beta$ is also a multi-index and $|\alpha-\beta|=|\alpha|-|\beta|$.

In this section, we will generalize those results obtained in Section 2 to the Hardy space $H^2(B_N)$ and the weighted Bergman spaces $A^2_s(B_N)$ ($s>-1$). In other words, we will compute the limits $$\lim\limits_{m\to\infty}||C_{\varphi}(e_m)||^2_{\mathcal{H}}\qquad \mbox{and} \qquad \lim\limits_{m\to\infty}||C^*_{\varphi}(e_m)||^2_{\mathcal{H}},$$ where the sequence $\{e_m\}$ is chosen to be  $$e_m(z)=\sqrt{\frac{\Gamma(t+m)}{\Gamma(t)!m!}}\biggl<z, \frac{a}{|a|}\biggr>^m$$
with $a=\varphi^{-1}(0)\ne 0$. Here, $t=N$ when $\mathcal{H}=H^2(B_N)$ and $t=N+s+1$ when $\mathcal{H}=A^2_s(B_N)$ ($s>-1$).
Surprisedly, when $\varphi$ is an involution automorphism of $B_N$, we find that the above two limits respectively have the same structures as those  on  the weighted Bergman space $A^2_s(D)$ ($s>-1$).

In fact, MacCluer and Pons \cite{MP06} have estimated the limits $$\lim\limits_{m\to\infty}||C_{\varphi}(z^m_2/||z^m_2||_{\mathcal{H}})||^2_{\mathcal{H}}\qquad \mbox{and} \qquad
\lim\limits_{m\to\infty}||C^*_{\varphi}(z^m_2/||z^m_2||_{\mathcal{H}})||^2_{\mathcal{H}}.$$ Comparing the two limits, they proved  that $C_\varphi$ is essentially normal on $\mathcal{H}$ if and only if $\varphi$ is unitary. Applying the method in the proof of Theorem 2.4, we can also compare the limits $\lim\limits_{m\to\infty}||C_{\varphi}(e_m)||^2_{\mathcal{H}}$ and $\lim\limits_{m\to\infty}||C^*_{\varphi}(e_m)||^2_{\mathcal{H}}$ to deduce the same result. Moreover, this will provide an important tool to investigate the essential normality of some linear fractional composition operators on  $\mathcal{H}$.

\begin{proposition}\label{prop 3.1}
Suppose that $\varphi_a(z)$ is an involution automorphism of $B_N$ interchanging $a$ and $0$. If $a\ne 0$, let
 $$e_m(z)=C_m\biggl<z, \frac{a}{|a|}\biggr>^m:=\sqrt{\frac{(N-1+m)!}{(N-1)!m!}}\biggl<z, \frac{a}{|a|}\biggr>^m,$$ then  on $H^2(B_N)$,
$$\lim\limits_{m\to\infty}||C_{\varphi_a}(e_m)||^2=(1-|a|^2)^N\sum\limits^{\infty}_{k=0} \biggl(\frac{\Gamma(N+k)}{k!\Gamma(N)}\biggr)^2|a|^{2k}.$$
\end{proposition}
\proof First, using the formula (see p.15 of \cite{Zhu})
$$\int_{\partial B_N}|<z, \zeta>|^{2m}d\sigma(\zeta)=\frac{(N-1)!m!}{(N-1+m)!}|z|^{2m},$$
we see that{\setlength\arraycolsep{2pt}
\begin{eqnarray*}||<z, a/|a|>^m||^2 &=&\int_{\partial B_N}|<\zeta, a/|a|>|^{2m}d\sigma(\zeta)
\\ &=&\frac{(N-1)!m!}{(N-1+m)!}\biggl|\frac{a}{|a|}\biggr|^{2m}=\frac{(N-1)!m!}{(N-1+m)!}.
\end{eqnarray*}}Thus $\{e_m\}$ is a  sequence of $B_N$ which weakly converges to zero with $||e_m||=1$.

Now, we have{\setlength\arraycolsep{2pt}
\begin{eqnarray*} ||C_{\varphi_a}(e_m)||^2 &= & C_m^2||<\varphi_a(z), a/|a|>^m||^2
\\ &=& C_m^2\int_{\partial B_N}|<\varphi_a(\zeta), a/|a|>|^{2m}d\sigma(\zeta)
\\ &=& C_m^2\int_{\partial B_N}|<\zeta, a/|a|>|^{2m}\biggl(\frac{1-|a|^2}{|1-<\zeta, a>|^2}\biggr)^Nd\sigma(\zeta)
\\ &=& (1-|a|^2)^N\frac{C_m^2}{|a|^{2m}}\int_{\partial B_N}|<\zeta, a>|^{2m}\biggl|\sum\limits^{\infty}_{k=0} \frac{\Gamma(N+k)}{k!\Gamma(N)}<\zeta, a>^k\biggr|^2d\sigma(\zeta)
\\ &=& (1-|a|^2)^N\frac{C_m^2}{|a|^{2m}}\int_{\partial B_N}\biggl|\sum\limits^{\infty}_{k=0} \frac{\Gamma(N+k)}{k!\Gamma(N)}<\zeta, a>^{k+m}\biggr|^2d\sigma(\zeta)
\\ &=& (1-|a|^2)^N\frac{C_m^2}{|a|^{2m}}\sum\limits^{\infty}_{k=0} \biggl(\frac{\Gamma(N+k)}{k!\Gamma(N)}\biggr)^2\int_{\partial B_N}|<\zeta, a>|^{2(k+m)}d\sigma(\zeta)
\\ &=& (1-|a|^2)^N\frac{C_m^2}{|a|^{2m}}\sum\limits^{\infty}_{k=0} \biggl(\frac{\Gamma(N+k)}{k!\Gamma(N)}\biggr)^2\frac{(N-1)!(m+k)!}{(N-1+m+k)!}|a|^{2(k+m)}
\\ &=& (1-|a|^2)^N\sum\limits^{\infty}_{k=0} \biggl(\frac{\Gamma(N+k)}{k!\Gamma(N)}\biggr)^2|a|^{2k}\frac{(N-1+m)!(m+k)!}{(N-1+m+k)!m!},
\end{eqnarray*}}where in the third line, we  have used the change of variables formula (see Corollary 4.4 in \cite{Zhu}), and in the sixth line, we  used the orthogonality of the functions $<z, a>^{k_1}$ and $<z, a>^{k_2}$   in $L^2(\partial B_N, d\sigma)$
when $k_1\ne k_2$.

Since $|a|<1$  and  for any fixed non-negative integer $k$,
$$\frac{(N-1+m)!(m+k)!}{(N-1+m+k)!m!}\to 1$$ as $m\to\infty$, we apply the dominated convergence theorem to get that{\setlength\arraycolsep{2pt}
\begin{eqnarray*} && \lim\limits_{m\to\infty}||C_{\varphi_a}(e_m)||^2
\\ &=& \lim\limits_{m\to\infty}(1-|a|^2)^N\sum\limits^{\infty}_{k=0} \biggl(\frac{\Gamma(N+k)}{k!\Gamma(N)}\biggr)^2|a|^{2k}\frac{(N-1+m)!(m+k)!}{(N-1+m+k)!m!}
\\ &=&(1-|a|^2)^N\sum\limits^{\infty}_{k=0} \biggl(\frac{\Gamma(N+k)}{k!\Gamma(N)}\biggr)^2|a|^{2k}.
\end{eqnarray*}}
 $\Box$
\\ \par Similar computation gives the following result for the weighted Bergman space $A^2_s(B_N)$ ($s>-1$). We only need use the change of variables formula for $A^2_s(B_N)$  (see Proposition 1.13 in \cite{Zhu}), so we omit its proof.
\begin{proposition}\label{prop 3.2}
Suppose that $\varphi_a(z)$ is an involution automorphism of $B_N$ interchanging $a$ and $0$. If $a\ne 0$, let
 $$e_m(z)=\sqrt{\frac{\Gamma(N+s+1+m)}{\Gamma(N+s+1)!m!}}\biggl<z, \frac{a}{|a|}\biggr>^m,$$ then  on $A^2_s(B_N)$ ($s>-1$),
$$\lim\limits_{m\to\infty}||C_{\varphi_a}(e_m)||^2_s=(1-|a|^2)^{N+s+1}\sum\limits^{\infty}_{k=0} \biggl(\frac{\Gamma(N+s+1+k)}{k!\Gamma(N+s+1)}\biggr)^2|a|^{2k}.$$
\end{proposition}

\begin{lemma}\label{lem 3.3} For any positive integers $m$ and $k$, there exist positive numbers $a_1, a_2,\ldots, a_k$ such that
$$\frac{(m+k)!}{m!}=M_k+a_1M_{k-1}+\cdots+a_{k-1}M_1+a_k,\eqno(3.1)$$
where $M_i=m(m-1)\cdots[m-(i-1)]$, ($i=1,\ldots, k$).
\end{lemma}
\proof We use induction, when $k=2$, it is clear that $$\frac{(m+2)!}{m!}=M_2+a_1M_1+a_2,$$
where $M_2=m(m-1)$, $M_1=m$, $a_1=2^2=4$, $a_2=2!=2$.

Assume that (3.1) is true for $k$, we prove that it is also true for $k+1$. By induction assumption (write $a_0=1$ and $M_0=1$),
{\setlength\arraycolsep{2pt}
\begin{eqnarray*} \frac{(m+k+1)!}{m!}&=&\frac{(m+k)!}{m!}(m+k+1)
\\ &=&\biggl(\sum\limits^k_{i=0}a_iM_{k-i}\biggr)(m+k+1)
\\ &=&\sum\limits^k_{i=0}a_iM_{k-i}[m-(k-i)+(2k-i+1)]
\\ &=&\sum\limits^k_{i=0}a_iM_{k-i}[m-(k-i)]+\sum\limits^k_{i=0}(2k-i+1)a_iM_{k-i}
\\ &=&\sum\limits^k_{i=0}a_iM_{k-i+1}+\sum\limits^{k+1}_{i=1}(2k-i+2)a_{i-1}M_{k-i+1}
\\ &=& a_0M_{k+1}+\sum\limits^k_{i=1}[a_i+(2k-i+2)a_{i-1}]M_{k-i+1}+(k+1)a_kM_0
\\ &=& \sum\limits^{k+1}_{i=0}a_i'M_{(k+1)-i},
\end{eqnarray*}}where $a'_0=1$, $a'_i=a_i+(2k-i+2)a_{i-1}$, $a'_{k+1}=(k+1)a_k=(k+1)!$. Hence, the result is true for $k+1$. The proof is complete.
\ \  $\Box$

\begin{proposition}\label{prop 3.3}
Suppose that   $\varphi$ is an  automorphism of $B_N$ with  $a=\varphi^{-1}(0)\ne 0$. Let
 $$e_m(z)=C_m\biggl<z, \frac{a}{|a|}\biggr>^m:=\sqrt{\frac{\Gamma(t+m)}{\Gamma(t)!m!}}\biggl<z, \frac{a}{|a|}\biggr>^m,$$ then on the space $\mathcal{H}$,
$$\lim\limits_{m\to\infty}||C_\varphi(e_m)||_{\mathcal{H}}^2\le (1-|a|^2)^t\sum\limits^{\infty}_{k=0} \biggl(\frac{\Gamma(t+k)}{k!\Gamma(t)}\biggr)^2|a|^{2k},$$
where $t=N$ when $\mathcal{H}=H^2(B_N)$ and $t=N+s+1$ when $\mathcal{H}=A^2_s(B_N)$ ($s>-1$).
\end{proposition}
\proof In the proof, we only discuss the case  $\mathcal{H}=H^2(B_N)$. If $\varphi$ is an  automorphism of $B_N$ with  $\varphi(a)=0$, by Theorem 2.25 in \cite{Rudin} or Theorem 1.4 in \cite{Zhu}, we have $\varphi=U\varphi_a$, where $U$ is unitary and $\varphi_a$ is an involution automorphism of $B_N$ that interchanges $a$ and $0$. Notice that{\setlength\arraycolsep{2pt}
\begin{eqnarray*} ||C_\varphi(e_m)||^2 &=& C^2_m||<\varphi(z), a/|a|>^m||^2=C^2_m||<U\varphi_a(z), a/|a|>^m||^2
\\ &=& C^2_m||<\varphi_a(z), U^{-1}(a/|a|)>^m||^2
\\ &=& C^2_m||C_{\varphi_a}(<z, U^{-1}(a/|a|)>^m)||^2.
\end{eqnarray*}}Hence, if we can show that
$$ C^2_m||C_{\varphi_a}(<z, \eta>^m)||^2\le  C^2_m||C_{\varphi_a}(<z, a/|a|>^m)||^2=||C_{\varphi_a}(e_m)||^2$$
for any $\eta\in \partial B_N$, then the desired result follows from Proposition 3.1.

Now, for  $\eta\in \partial B_N$, we compute that{\setlength\arraycolsep{2pt}
\begin{eqnarray*} && ||C_{\varphi_a}(<z, \eta>^m)||^2
=  ||<\varphi_a(z), \eta>^m||^2
\\ &=& \int_{\partial B_N}|<\varphi_a(\zeta), \eta>|^{2m}d\sigma(\zeta)
\\ &=& \int_{\partial B_N}|<\zeta, \eta>|^{2m}\biggl(\frac{1-|a|^2}{|1-<\zeta, a>|^2}\biggr)^Nd\sigma(\zeta)
\\ &=& (1-|a|^2)^N\int_{\partial B_N}|<\zeta, \eta>|^{2m}\biggl|\sum\limits^{\infty}_{k=0} \frac{\Gamma(N+k)}{k!\Gamma(N)}<\zeta, a>^k\biggr|^2d\sigma(\zeta)
\\ &=& (1-|a|^2)^N\int_{\partial B_N}\biggl|\sum\limits^{\infty}_{k=0} \frac{\Gamma(N+k)}{k!\Gamma(N)}<\zeta, a>^k<\zeta,\eta>^m\biggr|^2d\sigma(\zeta)
\\ &=& (1-|a|^2)^N\sum\limits^{\infty}_{k=0} \biggl(\frac{\Gamma(N+k)}{k!\Gamma(N)}\biggr)^2\int_{\partial B_N}|<\zeta, a>^k<\zeta,\eta>^m|^2d\sigma(\zeta).
\end{eqnarray*}}In the above calculations, we have used the change of variables formula on $H^2(B_N)$ and the orthogonality of $<\zeta, a>^{k_1}<\zeta,\eta>^m$ and $<\zeta, a>^{k_2}<\zeta,\eta>^m$ when $k_1\ne k_2$ (It is easy to check that).

In order to better understand our technique, we first estimate the above integral when $k=1$, {\setlength\arraycolsep{2pt}
\begin{eqnarray*} && \int_{\partial B_N}|<\zeta, a><\zeta,\eta>^m|^2d\sigma(\zeta)=\int_{\partial B_N}\biggl|\sum\limits_{i=1}^N\overline{a_i}\zeta_i\sum\limits_{|\alpha|=m}\frac{m!}{\alpha!}\overline{\eta}^\alpha \zeta^\alpha\biggr|^2d\sigma(\zeta)
\\&=&  \int_{\partial B_N}\biggl<\sum\limits_{i=1}^N\overline{a_i}\zeta_i\sum\limits_{|\alpha|=m}\frac{m!}{\alpha!}\overline{\eta}^\alpha \zeta^\alpha, \sum\limits_{j=1}^N\overline{a_j}\zeta_j\sum\limits_{|\beta|=m}\frac{m!}{\beta!}\overline{\eta}^\beta \zeta^\beta\biggr>d\sigma(\zeta)
\\&=&  \int_{\partial B_N}\biggl(\sum\limits_{i=1}^N\sum\limits_{|\alpha|=m}\frac{m!}{\alpha!}\frac{m!}{\alpha!}|a_i|^2|\eta^\alpha|^2 |\zeta_i|^2 |\zeta^\alpha|^2+ \sum\limits_{i\ne j}\sum\limits_{\alpha\ne \beta}\frac{m!}{\alpha!}\frac{m!}{\beta!}\overline{a_i} a_j \overline{\eta}^\alpha \eta^\beta \zeta_i \overline{\zeta_j} \zeta^\alpha  \overline{\zeta}^\beta\biggr)d\sigma(\zeta)
\\ &=& \sum\limits_{i=1}^N\sum\limits_{|\alpha|=m}\frac{m!}{\alpha!}\frac{m!}{\alpha!}|a_i|^2|\eta^\alpha|^2 \frac{(N-1)!\alpha_1!\cdots(\alpha_i+1)!\cdots \alpha_N!}{(N-1+m+1)!}
\\ && + \sum\limits_{i\ne j}\sum\limits_{\substack{\alpha_i+1= \beta_i \\ \alpha_j= \beta_j+1 \\\alpha_l= \beta_l, l\ne i, j}}\frac{m!}{\alpha!}\frac{m!}{\beta!}\overline{a_i} a_j \overline{\eta}^\alpha \eta^\beta \frac{(N-1)!\beta_1!\cdots(\beta_j+1)!\cdots \beta_N!}{(N-1+m+1)!}
\\&=& \frac{(N-1)!m!}{(N-1+m+1)!}\biggl(\sum\limits_{i=1}^N\sum\limits_{|\alpha|=m}\frac{m!(\alpha_i+1)}{\alpha!}|a_i|^2|\eta^\alpha|^2
\\ && + \sum\limits_{i\ne j}\sum\limits_{\substack{\alpha_j\ge 1 \\ |\alpha|=m}}\frac{m!\alpha_j}{\alpha!}|\eta_1|^{2\alpha_1}\cdots|\eta_j|^{2(\alpha_j-1)}\cdots|\eta_N|^{2\alpha_N}  \overline{a_i} a_j\eta_i \overline{\eta_j} \biggr)
\end{eqnarray*}}Note that for any $i, j$,
{\setlength\arraycolsep{2pt}
\begin{eqnarray*} && |a_i|^2|\eta_j|^2+ |a_j|^2|\eta_i|^2
\ge 2 |a_i||a_j||\eta_i||\eta_j|
\\ && \ge 2\mbox{Re}\,(\overline{a_i} a_j\eta_i \overline{\eta_j})
 =\overline{a_i} a_j\eta_i \overline{\eta_j}+ a_i \overline{a_j}\, \overline{\eta_i}\eta_j.
\end{eqnarray*}}This gives{\setlength\arraycolsep{2pt}
\begin{eqnarray*} &&\sum\limits_{i\ne j}\sum\limits_{\substack{\alpha_j\ge 1 \\ |\alpha|=m}}\frac{m!\alpha_j}{\alpha!}\biggl(|a_i|^2|\eta^\alpha|^2
 -|\eta_1|^{2\alpha_1}\cdots|\eta_j|^{2(\alpha_j-1)}\cdots|\eta_N|^{2\alpha_N}  \overline{a_i} a_j\eta_i \overline{\eta_j}\biggr)
\\ & =& \sum\limits_{i\ne j}\sum\limits_{\substack{\alpha_j\ge 1 \\ |\alpha|=m}}\frac{m!\alpha_j}{\alpha!}|\eta_1|^{2\alpha_1}\cdots|\eta_j|^{2(\alpha_j-1)}\cdots|\eta_N|^{2\alpha_N}|a_i|^2|\eta_j|^2
\\ && -\sum\limits_{i\ne j}\sum\limits_{\substack{\alpha_j\ge 1 \\ |\alpha|=m}}\frac{m!\alpha_j}{\alpha!}|\eta_1|^{2\alpha_1}\cdots|\eta_j|^{2(\alpha_j-1)}\cdots|\eta_N|^{2\alpha_N}  \overline{a_i} a_j\eta_i \overline{\eta_j} \ge 0.
\end{eqnarray*}}On the other hand,
{\setlength\arraycolsep{2pt}
\begin{eqnarray*} &&\sum\limits_{i=1}^N\sum\limits_{|\alpha|=m}\frac{m!(\alpha_i+1)}{\alpha!}|a_i|^2|\eta^\alpha|^2
 + \sum\limits_{i\ne j}\sum\limits_{\substack{\alpha_j\ge 1 \\ |\alpha|=m}}\frac{m!\alpha_j}{\alpha!}|a_i|^2|\eta^\alpha|^2
 \\&=& \sum\limits_{i=1}^N\sum\limits_{|\alpha|=m}\frac{m!\alpha_i}{\alpha!}|a_i|^2|\eta^\alpha|^2
 +\sum\limits_{i=1}^N\sum\limits_{|\alpha|=m}\frac{m!}{\alpha!}|a_i|^2|\eta^\alpha|^2
 + \sum\limits_{i\ne j}\sum\limits_{\substack{\alpha_j\ge 1 \\ |\alpha|=m}}\frac{m!\alpha_j}{\alpha!}|a_i|^2|\eta^\alpha|^2
\\&=& m\sum\limits_{i=1}^N\sum\limits_{\substack{\alpha_i\ge 1 \\ |\alpha|=m}}\frac{(m-1)!}{\alpha_1!\cdots(\alpha_i-1)!\cdots \alpha_N!}|\eta_1|^{2\alpha_1}\cdots|\eta_i|^{2(\alpha_i-1)}\cdots|\eta_N|^{2\alpha_N}|a_i|^2|\eta_i|^2
 \\ && +m \sum\limits_{i\ne j}\sum\limits_{\substack{\alpha_j\ge 1 \\ |\alpha|=m}}\frac{(m-1)!}{\alpha_1!\cdots(\alpha_j-1)!\cdots \alpha_N!}|\eta_1|^{2\alpha_1}\cdots|\eta_j|^{2(\alpha_j-1)}\cdots|\eta_N|^{2\alpha_N}|a_i|^2|\eta_j|^2
\\ &&  +\sum\limits_{i=1}^N|a_i|^2\sum\limits_{|\alpha|=m}\frac{m!}{\alpha!}|\eta^\alpha|^2
\\ &=& m|\eta|^{2(m-1)}\biggl(\sum\limits_{i=1}^N|a_i|^2|\eta_i|^2+\sum\limits_{i\ne j}|a_i|^2|\eta_j|^2\biggr)+|a|^2|\eta|^{2m}
\\ &=& m |\eta|^{2(m-1)} |a|^2|\eta|^2+|a|^2|\eta|^{2m}
=(m+1)|a|^2|\eta|^{2m}.
\end{eqnarray*}}Therefore, for $\eta\in\partial B_N$ and $k=1$,
{\setlength\arraycolsep{2pt}
\begin{eqnarray*}&& \int_{\partial B_N}|<\zeta, a><\zeta,\eta>^m|^2d\sigma(\zeta)
\\ &=& \frac{(N-1)!m!}{(N-1+m+1)!}\biggl(\sum\limits_{i=1}^N\sum\limits_{|\alpha|=m}\frac{m!(\alpha_i+1)}{\alpha!}|a_i|^2|\eta^\alpha|^2
+ \sum\limits_{i\ne j}\sum\limits_{\substack{\alpha_j\ge 1 \\ |\alpha|=m}}\frac{m!\alpha_j}{\alpha!}|a_i|^2|\eta^\alpha|^2
\\ && -\sum\limits_{i\ne j}\sum\limits_{\substack{\alpha_j\ge 1 \\ |\alpha|=m}}\frac{m!\alpha_j}{\alpha!}|a_i|^2|\eta^\alpha|^2
 + \sum\limits_{i\ne j}\sum\limits_{\substack{\alpha_j\ge 1 \\ |\alpha|=m}}\frac{m!\alpha_j}{\alpha!}|\eta_1|^{2\alpha_1}\cdots|\eta_j|^{2(\alpha_j-1)}\cdots|\eta_N|^{2\alpha_N}  \overline{a_i} a_j\eta_i \overline{\eta_j} \biggr).
\\ &\le &\frac{(N-1)!m!}{(N-1+m+1)!}(m+1)|a|^2|\eta|^{2m}= \frac{(N-1)!(m+1)!}{(N-1+m+1)!}|a|^2.
\end{eqnarray*}}

Next, for any integer $k\ge 2$, applying similar arguments and more complicated calculations, we can deduce that{\setlength\arraycolsep{2pt}
\begin{eqnarray*}&& \int_{\partial B_N}|<\zeta, a>^k<\zeta,\eta>^m|^2d\sigma(\zeta)
\\ &=&\int_{\partial B_N}\biggl|\sum\limits_{|\gamma|=k}\frac{k!}{\gamma!}\overline{a}^\gamma \zeta^\gamma\sum\limits_{|\alpha|=m}\frac{m!}{\alpha!}\overline{\eta}^\alpha \zeta^\alpha\biggr|^2d\sigma(\zeta)
\\ &=&\int_{\partial B_N}\biggl<\sum\limits_{|\gamma|=k}\frac{k!}{\gamma!}\overline{a}^\gamma \zeta^\gamma\sum\limits_{|\alpha|=m}\frac{m!}{\alpha!}\overline{\eta}^\alpha \zeta^\alpha, \sum\limits_{|\delta|=k}\frac{k!}{\delta!}\overline{a}^\delta \zeta^\delta\sum\limits_{|\beta|=m}\frac{m!}{\beta!}\overline{\eta}^\beta \zeta^\beta\biggr>d\sigma(\zeta)
\\ &=&\int_{\partial B_N}\biggl(\sum\limits_{\gamma=\delta}\sum\limits_{\alpha=\beta}\frac{k!}{\gamma!}\frac{k!}{\gamma!}\frac{m!}{\alpha!}\frac{m!}{\alpha!}|a^\gamma|^2
|\eta^\alpha|^2| \zeta^{\gamma+\alpha}|^2
\\ && +\sum\limits_{\gamma\ne\delta}\sum\limits_{\alpha\ne\beta}\frac{k!}{\gamma!}\frac{k!}{\delta!}\frac{m!}{\alpha!}\frac{m!}{\beta!}\overline{a}^\gamma
a^\delta\overline{\eta}^\alpha \eta^\beta \zeta^{\gamma+\alpha} \overline{\zeta}^{\delta+\beta}
\biggr)d\sigma(\zeta)
\\ &=&\sum\limits_{\gamma=\delta}\sum\limits_{\alpha=\beta}\frac{k!}{\gamma!}\frac{k!}{\gamma!}\frac{m!}{\alpha!}\frac{m!}{\alpha!}|a^\gamma|^2
|\eta^\alpha|^2\frac{(N-1)!(\gamma+\alpha)!}{(N-1+m+k)!}
\\ &&+\sum\limits_{\gamma\ne\delta}\sum\limits_{\gamma+\alpha=\delta+\beta}\frac{k!}{\gamma!}\frac{k!}{\delta!}\frac{m!}{\alpha!}\frac{m!}{\beta!}
\overline{a}^\gamma a^\delta\overline{\eta}^\alpha \eta^\beta \frac{(N-1)!(\gamma+\alpha)!}{(N-1+m+k)!}
\\ &=& \frac{(N-1)!m!}{(N-1+m+k)!} \biggl( \sum\limits_{\gamma=\delta}\sum\limits_{|\alpha|=m}\frac{k!}{\gamma!}\frac{k!}{\gamma!}\frac{m!}{\alpha!}\frac{(\gamma+\alpha)!}{\alpha!}|a^\gamma|^2
|\eta^\alpha|^2
\\ &&+\sum\limits_{\gamma\ne\delta}\sum\limits_{\gamma+\alpha\ge \delta}\frac{k!}{\gamma!}\frac{k!}{\delta!}\frac{m!}{\alpha!}\frac{(\gamma+\alpha)!}{(\gamma+\alpha-\delta)!}
\overline{a}^\gamma a^\delta\overline{\eta}^\alpha \eta^{\gamma+\alpha-\delta}\biggr)
\\ &\le &\frac{(N-1)!m!}{(N-1+m+k)!}(m+1)\cdots(m+k)|a|^{2k}|\eta|^{2m}
\\ &=& \frac{(N-1)!(m+k)!}{(N-1+m+k)!}|a|^{2k}.
\end{eqnarray*}}Here, in order to handle  $$\frac{(\gamma+\alpha)!}{\alpha!}=\frac{(\alpha_1+\gamma_1)!}{\alpha_1!}\cdot\frac{(\alpha_2+\gamma_2)!}{\alpha_2!}\cdots\frac{(\alpha_N+\gamma_N)!}{\alpha_N!}$$ and
$$\frac{(\gamma+\alpha)!}{(\gamma+\alpha-\delta)!}=\frac{(\gamma_1+\alpha_1)!}{(\gamma_1+\alpha_1-\delta_1)!}\cdot
\frac{(\gamma_2+\alpha_2)!}{(\gamma_2+\alpha_2-\delta_2)!}\cdots\frac{(\gamma_N+\alpha_N)!}{(\gamma_N+\alpha_N-\delta_N)!},$$ we  have used Lemma 3.3. We omit the details of computations, which are too complicated to display here. However, we include a computation for the case $N=2$ and $k=2$ in the Appendix to illustrate some further details of the idea.

Finally, combining the above conclusion  with the proof of Proposition 3.1, we get{\setlength\arraycolsep{2pt}
\begin{eqnarray*} && C^2_m||C_{\varphi_a}(<z, \eta>^m)||^2
\\ &= & C^2_m (1-|a|^2)^N\sum\limits^{\infty}_{k=0} \biggl(\frac{\Gamma(N+k)}{k!\Gamma(N)}\biggr)^2\int_{\partial B_N}|<\zeta, a>^k<\zeta,\eta>^m|^2d\sigma(\zeta)
\\ & \le &(1-|a|^2)^N\sum\limits^{\infty}_{k=0} \biggl(\frac{\Gamma(N+k)}{k!\Gamma(N)}\biggr)^2 |a|^{2k}\frac{(N-1+m)!(m+k)!}{(N-1+m+k)!m!}
\\ &=&||C_{\varphi_a}(e_m)||^2.
\end{eqnarray*}}
$\Box$

\begin{proposition}\label{prop 3.4}
Suppose that   $\varphi$ is an  automorphism of $B_N$ with  $a=\varphi^{-1}(0)\ne 0$. Let
 $$e_m(z)=C_m\biggl<z, \frac{a}{|a|}\biggr>^m:=\sqrt{\frac{(N-1+m)!}{(N-1)!m!}}\biggl<z, \frac{a}{|a|}\biggr>^m.$$ Then on $H^2(B_N)$,
$$\lim\limits_{m\to\infty}||C_\varphi^*(e_m)||^2=\frac{1}{(1-|a|^2)^N}\sum\limits^{N}_{k=0}\left(\begin{array}{c}
N
\\ k
\end{array}\right)^2|a|^{2k}.$$
\end{proposition}
\proof As  in the proof of Proposition 3.4, we  may write  $\varphi=U\varphi_a$. Thus, $C_\varphi=C_{\varphi_a}C_U$ and $C_\varphi C^*_\varphi= C_{\varphi_a}C_UC^*_UC^*_{\varphi_a}=C_{\varphi_a}C^*_{\varphi_a}$. Applying Proposition 3.1, we have
$$C_\varphi C^*_\varphi=C_{\varphi_a}C^*_{\varphi_a}=T_f+K,$$  where $T_f$ is the Toeplitz operator with symbol
$$f(z)=\biggl(\frac{|1-<z,a>|^2}{1-|a|^2}\biggr)^N$$
 and $K$ is a compact operator on $H^2(B_N)$. Write{\setlength\arraycolsep{2pt}
\begin{eqnarray*} f(z)&=&\frac{1}{(1-|a|^2)^N}(1-<z,a>)^N(1-\overline{<z,a>})^N
\\ &=& \frac{1}{(1-|a|^2)^N}\sum\limits^{N}_{k=0}\left(\begin{array}{c}
N
\\ k
\end{array}\right)(-1)^k<z,a>^k\sum\limits^{N}_{l=0}\left(\begin{array}{c}
N
\\ l
\end{array}\right)(-1)^l\overline{<z,a>}^l
\\ &=& \frac{1}{(1-|a|^2)^N}\sum\limits^{N}_{k=0}\sum\limits^{N}_{l=0}\left(\begin{array}{c}
N
\\ k
\end{array}\right)\left(\begin{array}{c}
N
\\ l
\end{array}\right)(-1)^{k+l}<z,a>^k\overline{<z,a>}^l
\end{eqnarray*}}and{\setlength\arraycolsep{2pt}
\begin{eqnarray*} && <z,a>^k\overline{<z,a>}^l<z,a>^m=\sum\limits_{|\gamma|=k}\frac{k!}{\gamma!}\overline{a}^\gamma z^\gamma \sum\limits_{|\delta|=l}\frac{l!}{\delta!}a^\delta \overline{z}^\delta \sum\limits_{|\alpha|=m}\frac{m!}{\alpha!}\overline{a}^\alpha z^\alpha
\\ &=& \sum\limits_{|\gamma|=k}\sum\limits_{|\delta|=l}\sum\limits_{|\alpha|=m}\frac{k!}{\gamma!}\frac{l!}{\delta!}\frac{m!}{\alpha!}
\overline{a}^{\gamma+\alpha}a^\delta z^{\gamma+\alpha}\overline{z}^\delta.
\end{eqnarray*}}Note that  (see Equation (12) in \cite{MP06})
$$P(\overline{z}^\delta z^{\gamma+\alpha})=T^*_{z^\delta}(z^{\gamma+\alpha})=\frac{(N+|\gamma+\alpha|-|\delta|-1)!}{(N+|\gamma+\alpha|-1)!}
\cdot\frac{(\gamma+\alpha)!}{(\gamma+\alpha-\delta)!}z^{\gamma+\alpha-\delta}$$
when $\gamma_i+\alpha_i\ge \delta_i$ for $1\le i\le N$ and $0$ otherwise, where $P$ denotes the orthogonal projection of $L^2(\partial B_N, d\sigma)$ onto $H^2(B_N)$. Combining $|\gamma+\alpha-\delta|=|\gamma|+|\alpha|-|\delta|=k+m-l$ with the orthogonality of the functions $<z, a>^{k_1}$ and $<z, a>^{k_2}$   in $L^2(\partial B_N, d\sigma)$
when $k_1\ne k_2$, we see that the inner product $$P[ <z,a>^k\overline{<z,a>}^l<z,a>^m]$$ with $<z,a>^m$ is zero unless $k=l$. This yields that{\setlength\arraycolsep{2pt}
\begin{eqnarray*} && <T_f(e_m), e_m> = \frac{C_m^2}{|a|^{2m}}<T_f(<z,a>^m), <z,a>^m>
\\ &=& \frac{1}{(1-|a|^2)^N}\frac{C_m^2}{|a|^{2m}}\sum\limits^{N}_{k=0}\left(\begin{array}{c}
N
\\ k
\end{array}\right)^2 <P[ <z,a>^k\overline{<z,a>}^k<z,a>^m], <z,a>^m>.
\end{eqnarray*}}

First, when $k=1$, we compute the projection{\setlength\arraycolsep{2pt}
\begin{eqnarray*} &&P[ <z,a>\overline{<z,a>}<z,a>^m]
\\ &=& \sum\limits_{i=1}^N\sum\limits_{j=1}^N\sum\limits_{|\alpha|=m}\frac{m!}{\alpha!}
\overline{a_i} a_j \overline{a}^\alpha P(\overline{z_j} z_i z^\alpha)
\\ &=& \sum\limits_{i=1}^N\sum\limits_{|\alpha|=m}\frac{m!}{\alpha!}
\overline{a}^\alpha |a_i|^2 \frac{\alpha_i+1}{N+m}z^\alpha+\sum\limits_{i\ne j}\sum\limits_{\substack{\alpha_j\ge 1 \\ |\alpha|=m}}\frac{m!}{\alpha!}
\overline{a}^\alpha \overline{a_i} a_j  \frac{\alpha_j}{N+m}\frac{ z_i z^\alpha}{z_j}
\end{eqnarray*}}and  the inner product{\setlength\arraycolsep{2pt}
\begin{eqnarray*} &&<P[ <z,a>\overline{<z,a>}<z,a>^m], <z,a>^m>
\\ &=& \biggl<\sum\limits_{i=1}^N\sum\limits_{|\alpha|=m}\frac{m!}{\alpha!}
\overline{a}^\alpha |a_i|^2 \frac{\alpha_i+1}{N+m}z^\alpha, \sum\limits_{|\beta|=m}\frac{m!}{\beta!}\overline{a}^\beta z^\beta\biggr>
\\ && +\sum\limits_{i\ne j}\sum\limits_{\substack{\alpha_j\ge 1 \\ |\alpha|=m}}\frac{m!}{\alpha!}
\overline{a}^\alpha \overline{a_i} a_j  \frac{\alpha_j}{N+m}\frac{ z_i z^\alpha}{z_j},
 \sum\limits_{|\beta|=m}\frac{m!}{\beta!}\overline{a}^\beta z^\beta\biggr>
\\ &=& \sum\limits_{i=1}^N\sum\limits_{|\alpha|=m}\frac{m!}{\alpha!}\frac{m!}{\alpha!}
|a^\alpha|^2 |a_i|^2 \frac{\alpha_i+1}{N+m}||z^\alpha||^2
\\ &&+\sum\limits_{i\ne j}\sum\limits_{\substack{\alpha_i+1=\beta_i \\ \alpha_j-1=\beta_j\\ \alpha_l=\beta_l, l\ne i,j}}\frac{m!}{\alpha!}
\frac{m!}{\beta!}\overline{a_i} a_j \overline{a}^\alpha  a^\beta \frac{\alpha_j}{N+m} ||z^\beta||^2
\\ &=& \sum\limits_{i=1}^N\sum\limits_{|\alpha|=m}\frac{m!}{\alpha!}\frac{m!}{\alpha!}
|a^\alpha|^2 |a_i|^2\frac{\alpha_i+1}{N+m} \cdot \frac{(N-1)!\alpha!}{(N-1+m)!}
\\ &&+\sum\limits_{i\ne j}\sum\limits_{\substack{\alpha_i+1=\beta_i \\ \alpha_j-1=\beta_j\\ \alpha_l=\beta_l, l\ne i,j}}\frac{m!}{\alpha!}
\frac{m!}{\beta!}\overline{a_i} a_j \overline{a}^\alpha  a^\beta \frac{\alpha_j}{N+m}\cdot \frac{(N-1)!\beta!}{(N-1+m)!}
\\ &=& \frac{(N-1)!m!}{(N-1+m)!}\biggl(\sum\limits_{i=1}^N\sum\limits_{|\alpha|=m}\frac{m!}{\alpha!}|a^\alpha|^2 |a_i|^2 \frac{\alpha_i}{N+m}
 +\sum\limits_{i=1}^N\sum\limits_{|\alpha|=m}\frac{m!}{\alpha!}|a^\alpha|^2 |a_i|^2 \frac{1}{N+m}
\end{eqnarray*}}
{\setlength\arraycolsep{2pt}
\begin{eqnarray*}&& +\sum\limits_{i\ne j}\sum\limits_{\substack{\alpha_j\ge 1\\ |\alpha|=m}}\frac{m!}{\alpha!}
 |a^\alpha|^2 |a_i|^2  \frac{\alpha_j}{N+m}\biggr)
 \\ &=& \frac{1}{C_m^2(N+m)}\biggl( m\sum\limits_{i=1}^N\sum\limits_{\substack{\alpha_i\ge 1 \\ |\alpha|=m}}\frac{(m-1)!}{\alpha_1!\cdots(\alpha_i-1)!\cdots \alpha_N!}|a_1|^{2\alpha_1}\cdots|a_i|^{2(\alpha_i-1)}\cdots|a_N|^{2\alpha_N}|a_i|^4
 \\ && +\sum\limits_{i=1}^N|a_i|^2\sum\limits_{|\alpha|=m}\frac{m!}{\alpha!}|a^\alpha|^2
 \\ && + m\sum\limits_{i\ne j}\sum\limits_{\substack{\alpha_j\ge 1 \\ |\alpha|=m}}\frac{(m-1)!}{\alpha_1!\cdots(\alpha_j-1)!\cdots \alpha_N!}|a_1|^{2\alpha_1}\cdots|a_j|^{2(\alpha_j-1)}\cdots|a_N|^{2\alpha_N}|a_i|^2|a_j|^2
 \biggr)
 \\ &=& \frac{1}{C_m^2(N+m)}\biggl[ m|a|^{2(m-1)}\biggl(\sum\limits_{i=1}^N|a_i|^2|a_i|^2+\sum\limits_{i\ne j}|a_i|^2|a_j|^2\biggr)
+  |a|^{2(m+1)} \biggr]
 \\ &=& \frac{1}{C_m^2(N+m)}(m+1)|a|^{2(m+1)}=\frac{|a|^{2m}}{C^2_m}\cdot\frac{m+1}{N+m}|a|^2,
\end{eqnarray*}}where in the above calculation, we have used the norm
$$||z^\alpha||^2=\frac{(N-1)!\alpha!}{(N-1+|\alpha|)!}.$$

Now, when  $k\ge 2$, using Lemma 3.3 and  similar idea in the proof of Proposition 3.4, we get
that{\setlength\arraycolsep{2pt}
\begin{eqnarray*} && P[ <z,a>^k\overline{<z,a>}^k<z,a>^m]
\\ &=& \sum\limits_{|\gamma|=k}\sum\limits_{|\delta|=k}\sum\limits_{|\alpha|=m}\frac{k!}{\gamma!}\frac{k!}{\delta!}\frac{m!}{\alpha!}
\overline{a}^{\gamma+\alpha}a^\delta P(\overline{z}^\delta z^{\gamma+\alpha})
\\ &=& \sum\limits_{|\gamma|=k}\sum\limits_{|\alpha|=m}\frac{k!}{\gamma!}\frac{k!}{\gamma!}\frac{m!}{\alpha!}
|a^\gamma|^2\overline{a}^\alpha \frac{(N+m-1)!}{(N+m+k-1)!}\cdot\frac{(\gamma+\alpha)!}{\alpha!}z^\alpha
\\ && + \sum\limits_{\substack{\gamma\ne \delta\\ |\gamma|=|\delta|=k}}\sum\limits_{\substack{\gamma+\alpha\ge \delta \\ |\alpha|=m}}\frac{k!}{\gamma!}\frac{k!}{\delta!}
\frac{m!}{\alpha!}\overline{a}^{\gamma+\alpha}a^\delta \frac{(N+m-1)!}{(N+m+k-1)!}\cdot\frac{(\gamma+\alpha)!}{(\gamma+\alpha-\delta)!}z^{\gamma+\alpha-\delta}
\end{eqnarray*}}
and{\setlength\arraycolsep{2pt}
\begin{eqnarray*} &&<P[ <z,a>^k\overline{<z,a>}^k<z,a>^m], <z,a>^m>
\\ &=& \biggl<P[ <z,a>^k\overline{<z,a>}^k<z,a>^m], \sum\limits_{|\beta|=m}\frac{m!}{\beta!}\overline{a}^\beta z^\beta\biggr>
\\ &=& \sum\limits_{ |\gamma|=k}\sum\limits_{|\alpha|=m}\frac{k!}{\gamma!}\frac{k!}{\gamma!}\frac{m!}{\alpha!}\frac{m!}{\alpha!}|a^\gamma|^2|a^\alpha|^2 \frac{(N+m-1)!}{(N+m+k-1)!}\cdot\frac{(\gamma+\alpha)!}{\alpha!}||z^\alpha||^2
\end{eqnarray*}}
{\setlength\arraycolsep{2pt}
\begin{eqnarray*} && + \sum\limits_{\substack{\gamma\ne \delta\\ |\gamma|=|\delta|=k}}\sum\limits_{\substack{\gamma+\alpha= \delta+\beta \\ |\alpha|=|\beta|=m}}\frac{k!}{\gamma!}\frac{k!}{\delta!}
\frac{m!}{\alpha!}\frac{m!}{\beta!}\overline{a}^{\gamma+\alpha}a^{\delta+\beta} \frac{(N+m-1)!}{(N+m+k-1)!}\cdot\frac{(\gamma+\alpha)!}{(\gamma+\alpha-\delta))!}||z^\beta||^2
\\ &=& \frac{1}{C_m^2}\frac{(N+m-1)!}{(N+m+k-1)!}\biggl(\sum\limits_{ |\gamma|=k}\sum\limits_{|\alpha|=m}\frac{k!}{\gamma!}\frac{k!}{\gamma!}\frac{m!}{\alpha!}\frac{(\gamma+\alpha)!}{\alpha!}|a^\gamma|^2|a^\alpha|^2
\\ && +\sum\limits_{\substack{\gamma\ne \delta\\ |\gamma|=|\delta|=k}}\sum\limits_{\substack{\gamma+\alpha\ge \delta\\ |\alpha|=m}}\frac{k!}{\gamma!}\frac{k!}{\delta!}\frac{m!}{\alpha!}\frac{(\gamma+\alpha)!}{(\gamma+\alpha-\delta))!}|a^{\gamma+\alpha}|^2\biggr)
 \\ &=&  \frac{1}{C_m^2}\frac{(N+m-1)!}{(N+m+k-1)!}(m+1)\cdots(m+k)|a|^{2(m+k)}
 \\ &= & \frac{|a|^{2m}}{C_m^2} \frac{(N+m-1)!(m+k)!}{(N+m+k-1)!m!}|a|^{2k}.
\end{eqnarray*}}Therefore, all above arguments give that{\setlength\arraycolsep{2pt}
\begin{eqnarray*} && <T_f(e_m), e_m>
\\ &=& \frac{1}{(1-|a|^2)^N}\frac{C_m^2}{|a|^{2m}}\sum\limits^{N}_{k=0}\left(\begin{array}{c}
N
\\ k
\end{array}\right)^2 <P[ <z,a>^k\overline{<z,a>}^k<z,a>^m], <z,a>^m>
\\ &=&  \frac{1}{(1-|a|^2)^N}\frac{C_m^2}{|a|^{2m}}\sum\limits^{N}_{k=0}\left(\begin{array}{c}
N
\\ k
\end{array}\right)^2\frac{|a|^{2m}}{C_m^2}  \frac{(N+m-1)!(m+k)!}{(N+m+k-1)!m!}|a|^{2k}
\\ &=&  \frac{1}{(1-|a|^2)^N}\sum\limits^{N}_{k=0}\left(\begin{array}{c}
N
\\ k
\end{array}\right)^2 \frac{(N+m-1)!(m+k)!}{(N+m+k-1)!m!}|a|^{2k}
.
\end{eqnarray*}}

At last, since $\{e_m\}$ is a sequence  which converges uniformly  to zero on any compact subset of $B_N$ and  $|a|<1$, we apply the dominated convergence theorem to see that{\setlength\arraycolsep{2pt}
\begin{eqnarray*} && \lim\limits_{m\to\infty}||C_\varphi^*(e_m)||^2= \lim\limits_{m\to\infty} <T_f(e_m), e_m>
\\ &=&   \lim\limits_{m\to\infty}\frac{1}{(1-|a|^2)^N}\sum\limits^{N}_{k=0}\left(\begin{array}{c}
N
\\ k
\end{array}\right)^2 \frac{(N+m-1)!(m+k)!}{(N+m+k-1)!m!}|a|^{2k}
\\ &=&   \frac{1}{(1-|a|^2)^N}\sum\limits^{N}_{k=0}\left(\begin{array}{c}
N
\\ k
\end{array}\right)^2 |a|^{2k}
.
\end{eqnarray*}}
$\Box$

\begin{proposition}\label{prop 3.5}
Suppose that   $\varphi$ is an  automorphism of $B_N$ with  $a=\varphi^{-1}(0)\ne 0$. Let
 $$e_m(z)=C_m\biggl<z, \frac{a}{|a|}\biggr>^m:=\sqrt{\frac{\Gamma(N+s+1+m)}{\Gamma(N+s+1)m!}}\biggl<z, \frac{a}{|a|}\biggr>^m.$$ Then on $A_s^2(B_N)$ ($s>-1$),
$$\lim\limits_{m\to\infty}||C_\varphi^*(e_m)||_s^2=\frac{1}{(1-|a|^2)^{N+s+1}}\sum\limits^{\infty}_{k=0}\left(\begin{array}{c}
N+s+1
\\ k
\end{array}\right)^2|a|^{2k}.$$
\end{proposition}
\proof Using similar method as in the proof of Proposition 3.5, we can obtain the desired result. For completeness, we  give an outline for its proof. On
$A_s^2(B_N)$ ($s>-1$), Proposition 2.1 gives
{\setlength\arraycolsep{2pt}
\begin{eqnarray*} && f(z)= \biggl(\frac{|1-<z,a>|^2}{1-|a|^2}\biggr)^{N+s+1}
\\ &=&
\frac{1}{(1-|a|^2)^{N+s+1}}(1-<z,a>)^{N+s+1}(1-\overline{<z,a>})^{N+s+1}
\\ &=& \frac{1}{(1-|a|^2)^{N+s+1}}\sum\limits^{\infty}_{k=0}\sum\limits^{\infty}_{l=0}\left(\begin{array}{c}
N+s+1
\\ k
\end{array}\right)\left(\begin{array}{c}
N+s+1
\\ l
\end{array}\right)(-1)^{k+l}<z,a>^k\overline{<z,a>}^l.
\end{eqnarray*}}An easy computation yields that
$$T^*_{z^\delta}(z^{\gamma+\alpha})=\frac{\Gamma(N+s+1+|\gamma+\alpha|-|\delta|)}{\Gamma(N+s+1+|\gamma+\alpha|)}
\cdot\frac{(\gamma+\alpha)!}{(\gamma+\alpha-\delta)!}z^{\gamma+\alpha-\delta}$$
when $\gamma_i+\alpha_i\ge \delta_i$ for $1\le i\le N$ and $0$ otherwise. Moreover,
we can calculate that{\setlength\arraycolsep{2pt}
\begin{eqnarray*} &&<P[ <z,a>^k\overline{<z,a>}^k<z,a>^m], <z,a>^m>
 \\ &=&  \frac{1}{C_m^2}\frac{\Gamma(N+s+1+m)}{\Gamma(N+s+1+m+k)}(m+1)\cdots(m+k)|a|^{2(m+k)}
 \\ &= & \frac{|a|^{2m}}{C_m^2} \frac{\Gamma(N+s+1+m)(m+k)!}{\Gamma(N+s+1+m+k)m!}|a|^{2k},
\end{eqnarray*}}where $P$ denotes the orthogonal projection of $L^2(B_N, dv_s)$ onto $A^2_s(B_N)$. Therefore, {\setlength\arraycolsep{2pt}
\begin{eqnarray*} && <T_f(e_m), e_m>=\frac{1}{(1-|a|^2)^{N+s+1}}\frac{C_m^2}{|a|^{2m}}<T_f(<z,a>^m), <z,a>^m>
\\ &=& \frac{1}{(1-|a|^2)^{N+s+1}}\frac{C_m^2}{|a|^{2m}}\times
\\ && \sum\limits^{\infty}_{k=0}\left(\begin{array}{c}
N+s+1
\\ k
\end{array}\right)^2<P[ <z,a>^k\overline{<z,a>}^k<z,a>^m], <z,a>^m>
\\ &=& \frac{1}{(1-|a|^2)^{N+s+1}}\frac{C_m^2}{|a|^{2m}}
\sum\limits^{\infty}_{k=0}\left(\begin{array}{c}
N+s+1
\\ k
\end{array}\right)^2\frac{|a|^{2m}}{C_m^2} \frac{\Gamma(N+s+1+m)(m+k)!}{\Gamma(N+s+1+m+k)m!}|a|^{2k}
\\ &=& \frac{1}{(1-|a|^2)^{N+s+1}}\sum\limits^{\infty}_{k=0}\left(\begin{array}{c}
N+s+1
\\ k
\end{array}\right)^2 \frac{\Gamma(N+s+1+m)(m+k)!}{\Gamma(N+s+1+m+k)m!}|a|^{2k}
.
\end{eqnarray*}}and whence{\setlength\arraycolsep{2pt}
\begin{eqnarray*} && \lim\limits_{m\to\infty}||C_\varphi^*(e_m)||^2_s= \lim\limits_{m\to\infty} <T_f(e_m), e_m>
\\ &=&   \lim\limits_{m\to\infty} \frac{1}{(1-|a|^2)^{N+s+1}}\sum\limits^{\infty}_{k=0}\left(\begin{array}{c}
N+s+1
\\ k
\end{array}\right)^2 \frac{\Gamma(N+s+1+m)(m+k)!}{\Gamma(N+s+1+m+k)m!}|a|^{2k}
\\ &=&  \frac{1}{(1-|a|^2)^{N+s+1}}\sum\limits^{\infty}_{k=0}\left(\begin{array}{c}
N+s+1
\\ k
\end{array}\right)^2|a|^{2k}
.
\end{eqnarray*}}
$\Box$
\\ \par
Now, using similar idea as in the proof of Theorem 2.4 and Proposition 3.4-3.6, we can deduce the following result, the proof will be omitted.

\begin{theorem}\label{the 3.6} Let $\varphi$ be  an automorphism  of $B_N$ with  $a=\varphi^{-1}(0)\ne 0$ and
 $$e_m(z)=\sqrt{\frac{\Gamma(t+m)}{\Gamma(t)!m!}}\biggl<z, \frac{a}{|a|}\biggr>^m.$$
Then on the space $\mathcal{H}$,
$$\lim\limits_{m\to\infty}(||C^*_\varphi(e_m)||_{\mathcal{H}}^2-||C_\varphi(e_m)||_{\mathcal{H}}^2)>0.$$
Here $t=N$ when $\mathcal{H}=H^2(B_N)$ and $t=N+s+1$ when $\mathcal{H}=A^2_s(B_N)$ ($s>-1$).
\end{theorem}

Therefore, as a  result of this theorem, we get that $C_\varphi$ is essentially normal on  $\mathcal{H}$ if and only if $\varphi$ is unitary.

\section {Essential normality of linear fractional  composition operators}
 In this section, we discuss two classes of linear fractional self-maps of $B_N$. They induce composition operators which are not essentially  normal  on the Hardy space  $H^2(B_N)$ or the weighted Bergman space $A^2_s(B_N)$ ($s>-1$).

In \cite{MW05},  MacCluer and Weir proved that if  the linear fractional map $\varphi$ restricted to the slice $[\zeta]=\{\lambda\zeta: \lambda\in\mathbb{C}\  \mbox{and}\ |\lambda|<1\}$ is a non-rotation automorphism of $[\zeta]$, then $C_\varphi$ is not essentially normal on $H^2(B_N)$ or $A^2_s(B_N)$ for any positive integer $s$ (see Theorem 7 of \cite{MW05}).
Using Theorem 3.7, we can generalize it  to the following result.

\begin{theorem}\label{the 4.1} If  $\varphi$ is   a linear fractional self-map   of $B_N$ and the restriction of $\varphi$ to
$$B_k=\{(z_1,\ldots, z_N)\in B_N;\  z_i=0\ \mbox{for}\ i> k\}$$
 is a non-rotation automorphism of $B_k$,  then $C_\varphi$ is not essentially normal on $H^2(B_N)$ or $A^2_s(B_N)$ for any real $s>-1$.
\end{theorem}

First, we need some lemmas.

\begin{lemma}\label{lem 4.2} Let $\varphi$   be a linear fractional self-map of $B_N$. If $\varphi$ maps $B_k$ into itself and the restriction of $\varphi$ to
$B_k$ is an automorphism, then the first $k$ coordinate functions of $\varphi$ and its adjoint map $\sigma$ depend only on the variables $z_1,\ldots, z_k$.
\end{lemma}
\proof The idea is similar to Lemma 3 of \cite{MW05}, but we still give a detail proof, because we need some parts below for the proof of Theorem 4.1. Assume
 that $$\varphi(z)=\frac{Az+B}{<z, C>+1}$$ with $A=(a_{ij})_{N\times N}$, $B=(b_i)_{N\times 1}$ and $C=(c_i)_{N\times 1}$. Then the adjoint map $\sigma$ must have the form $$\sigma(z)=\frac{A^*z-C}{<z, -B>+1}.$$

Since $\varphi$ is an automorphism when restricted to $B_k$, by Lemma 1 of \cite{MW05}, $\varphi$ and $\sigma$ all map $\partial B_k$ onto $\partial B_k$.
This implies that the last $N-k$ coordinate functions of $\varphi$ and $\sigma$ map the points of $\partial B_k$ to $0$. Thus, for $j>k$, the $j$-th coordinate
functions
$$\varphi_j(z)=\frac{a_{j1}z_1+\cdots a_{jN}z_N+b_j}{\overline{c_1}z_1+\cdots+\overline{c_N}z_N+1}$$
and $$\sigma_j(z)=\frac{\overline{a_{1j}}z_1+\cdots \overline{a_{Nj}}z_N-c_j}{-\overline{b_1}z_1-\cdots-\overline{b_N}z_N+1}$$
map the points $\lambda e_i$ ($i=1,\ldots, k$) with $\lambda\in \mathbb{C}$ and $|\lambda|=1$ to $0$, where $e_i$ ($i=1,\ldots, N$) form the normalized orthogonal  basis of $\mathbb{C}^N$. Hence $$a_{j1}=\cdots=a_{jk}=b_j=0$$  and $$\overline{a_{1j}}=\cdots=\overline{a_{kj}}=c_j=0$$
for $j>k$. That is,  $A$ must have the form
\begin{displaymath}
A=\left(
\begin{array}{cc}
A_1 & 0
\\ 0 & A_2
\end{array}\right)
\end{displaymath}with $A_1=(a_{ij})_{k\times k}$ and $A_2=(a_{ij})_{(N-k)\times (N-k)}$ and $B=(b_1,\ldots,b_k,0')=(B_1,0')\in \mathbb{C}^k\times\mathbb{C}^{N-k}$ and $C=(c_1,\ldots,c_k,0')=(C_1,0')\in \mathbb{C}^k\times\mathbb{C}^{N-k}$. From all these facts we deduce  the desired conclusions for $\varphi$ and $\sigma$. \ \  $\Box$
\\ \\
{\bf Lemma A.} (Lemma 1.9 of \cite{Zhu}) {\it Suppose $f$ is a function on $\partial B_N$ that depends only on  $z_1,\ldots, z_k$, where $1\le k<N$. Then $f$
 can be regarded as defined on $B_k$ and $$\int_{\partial B_N}f d\sigma=\left(
\begin{array}{c}
N-1
\\ k
\end{array}\right)\int_{B_k} (1-|w|^2)^{N-k-1}f(w)dv(w),$$ where $B_k$ is the unit ball in $\mathbb{C}^k$ and $dv$ is the normalized volume measure on $B_k$.}
\\ \par
Similarly, we have the following result for the weighted Bergman space $A^2_s(B_N)$ ($s>-1$), see exercise 4.27  in \cite{Zhu}. For convenience, we give a simple proof.
We identify $\mathbb{C}^N$ with $\mathbb{R}^{2N}$ using the real and imaginary parts of a complex number, and denote the usual Lebesgue measure on $\mathbb{C}^N$ by $dV=dx_1dy_1\cdots dx_N dy_N$. Then $dV=\frac{\pi^N}{N!}dv$  (see p.13 in \cite{Zhu}),  where $dv$ is the normalized volume measure on $B_N$.

\begin{lemma}\label{lem 4.3} Suppose $f$ is a function on $ B_N$ that depends only on  $z_1,\ldots, z_k$, where $1\le k<N$. Then for $s>-1$,
 $$\int_{ B_N}f(z) dv_s(z)=\int_{B_k} f(w)dv_{N-k+s}(w).$$
\end{lemma}
\proof  For $z\in\mathbb{C}^N$, write $z=(w, u)\in \mathbb{C}^k\times\mathbb{C}^{N-k}$. Applying Fubini's theorem and integration in polar coordinates (see 1.4.3 of \cite{Rudin} or Lemma 1.8 of \cite{Zhu}),
{\setlength\arraycolsep{2pt}
\begin{eqnarray*} &&  \int_{B_N} f(z)dv_s(z) = c_s\frac{N!}{\pi^N}\int_{B_N} f(z)(1-|z|^2)^s dV(z)
\\ &=& c_s\frac{N!}{\pi^N}\int_{B_k} f(w)dV(w)\int_{(1-|w|^2)B_{N-k}}(1-|w|^2-|u|^2)^sdV(u)
\\ &=& c_s\frac{N!}{\pi^N}\frac{\pi^{N-k}}{(N-k)!}\frac{\pi^k}{k!} \int_{B_k}f(w)(1-|w|^2)^sdv(w)\int_{(1-|w|^2)B_{N-k}}\biggl(1-\frac{|u|^2}{1-|w|^2}\biggr)^sdv(u)
\\ &=& c_s\frac{N!}{k!(N-k)!} \int_{B_k}f(w)(1-|w|^2)^sdv(w)\times
\\ && 2(N-k)\int_0^{1-|w|^2} r^{2(N-k)-1}\biggl(1-\frac{r^2}{1-|w|^2}\biggr)^sdr\int_{\partial B_{N-k}} d\sigma
\end{eqnarray*}}
{\setlength\arraycolsep{2pt}
\begin{eqnarray*}  &=& \frac{\Gamma(N+s+1)}{N!\Gamma(s+1)}\frac{N!}{k!(N-k)!}\frac{\Gamma(N-k+1)\Gamma(s+1)}{\Gamma(N-k+s+1)} \int_{B_k}f(w)(1-|w|^2)^{s+N-k}dv(w)
\\ &=& \frac{\Gamma(N+s+1)}{k!\Gamma(N-k+s+1)}\int_{B_k}f(w)(1-|w|^2)^{s+N-k}dv(w)
\\ &=& \int_{B_k} f(w)dv_{N-k+s}(w).
\end{eqnarray*}}
$\Box$
\\ \par
Based on Theorem 3.7 and the idea of Theorem 7 in \cite{MW05}, we  give the following proof for Theorem 4.1.
\\ \par
{\it Proof of Theorem 4.1.} First, we focus on the case of $H^2(B_N)$. In terms of the proof of Lemma 4.2, we find that $\varphi$ and its adjoint map $\sigma$
have the forms $$\varphi(z)=\biggl(\frac{A_1w+B_1}{<w, C_1>+1}, \frac{A_2u}{<w, C_1>+1}\biggr)$$ and
$$\sigma(z)=\biggl(\frac{A^*_1w-C_1}{<w, -B_1>+1}, \frac{A^*_2u}{<w, -B_1>+1}\biggr)$$
for $z=(w, u)\in \mathbb{C}^k\times\mathbb{C}^{N-k}$. Write
$$\widetilde{\varphi}(w)=\frac{A_1w+B_1}{<w, C_1>+1}\qquad\mbox{and}\qquad \widetilde{\sigma}(w)\frac{A^*_1w-C_1}{<w, -B_1>+1}.$$
Thus, by hypothesis,   $\widetilde{\varphi}$ is a non-rotation automorphism of $B_k$ and $\widetilde{\sigma}$ is its adjoint map.

Let
$a=\widetilde{\varphi}^{-1}(0)$ and $f_m(z)=C_m<w, a/|a|>^m=e_m(w)$, where $$C_m=\sqrt{\frac{\Gamma(N+m)}{\Gamma(N)!m!}}=\sqrt{\frac{\Gamma(k+(N-k-1)+m)}{\Gamma(k+(N-k-1))!m!}}.$$
By Lemma A, we have $$||f_m||^2=||e_m||_{N-k-1}^2=1,$$ so that $\{f_m\}$ is a normalized sequence on $B_N$ which weakly converges to $0$.

It is easy to check that $$C_\varphi(f_m)(z)=f_m(\varphi(z))=C_m<\widetilde{\varphi}(w), a/|a|>^m=C_{\widetilde{\varphi}}(e_m)(w).$$
On the other hand, we find that $C^*_\varphi=T_gC_\sigma T^*_h$ and $C^*_{\widetilde{\varphi}}=T_{\widetilde{g}}C_{\widetilde{\sigma}} T^*_{\widetilde{h}}$, where $$g(z)=\frac{1}{({<w, -B_1>+1})^N}=\widetilde{g}(w)$$ and $$h(z)=(<w, C_1>+1)^N=\widetilde{h}(w).$$ Moreover, for $\alpha\ge \beta$,
$T^*_{w^\beta}(w^\alpha)$ is identical when Toeplitz operator $T$ acts on $H^2(B_N)$ or $A^2_{N-k-1}(B_k)$.
It follows that{\setlength\arraycolsep{2pt}
\begin{eqnarray*}  C^*_\varphi(f_m)(z)&=&C_m T_gC_\sigma T^*_h(<w, a/|a|>^m)
\\ &=& C_mT_{\widetilde{g}}C_{\widetilde{\sigma}} T^*_{\widetilde{h}}(<w, a/|a|>^m)
\\ &=& C^*_{\widetilde{\varphi}}(e_m)(w).
\end{eqnarray*}}Using Lemma A again, we obtain
$$||C_\varphi(f_m)||^2=||C_{\widetilde{\varphi}}(e_m)||_{N-k-1}^2\quad\mbox{and}\quad||C^*_\varphi(f_m)||^2
=||C^*_{\widetilde{\varphi}}(e_m)||_{N-k-1}^2.$$
Now, by  Theorem 3.7, we get {\setlength\arraycolsep{2pt}
\begin{eqnarray*} &&
\lim\limits_{m\to\infty}(||C^*_\varphi(f_m)||^2-||C_\varphi(f_m)||^2)\\ &=&\lim\limits_{m\to\infty}
(||C^*_{\widetilde{\varphi}}(e_m)||_{N-k-1}^2-||C_{\widetilde{\varphi}}(e_m)||_{N-k-1}^2)>0.
\end{eqnarray*}}Hence $[C^*_\varphi, C_\varphi]$ is not compact on $H^2(B_N)$ and  the desired result holds.

On $A_s^2(B_N)$ ($s>-1$), all arguments are the same to those on $H^2(B_N)$, we only need use Lemma 4.3 to see that
$$||C_\varphi(f_m)||_s^2=||C_{\widetilde{\varphi}}(e_m)||_{N-k+s}^2\quad\mbox{and}\quad||C^*_\varphi(f_m)||^2_s=||C^*_{\widetilde{\varphi}}(e_m)||_{N-k+s}^2,$$
where $$f_m(z)=\frac{\Gamma(N+s+1+m)}{\Gamma(N+s+1)m!}<w,a/|a|>^m=e_m(w).$$
$\Box$
\\ \\
{\bf Remark 1.} Suppose that $\varphi$ is a linear fractional self-map of $B_N$. By the Cayley transform $\sigma_C(z)=(e_1+z)/(1-z_1)$, $z\in B_N$, $\varphi$ is conjugated to a map $\Phi$ on
$$H_N=\{(w_1,w')\in \mathbb{C}\times\mathbb{C}^{N-1}: \ \mbox{Re}\, w_1>|w'|^2\}$$ with the form
$$\Phi(z,u,v,w)=(z+2<u, a>+b, u+a, Dv, Aw),$$ where $(z,u,v,w)\in\mathbb{C}\times\mathbb{C}\times\mathbb{C}^{k-2}\times\mathbb{C}^{N-k}$ and
$\mbox{Re}\, b=|a|^2$, $D$ is unitary, $A$ is a matrix with $||A||<1$. Then using Theorem 4.1, we see that $C_\varphi$ is not essentially normal on the space $\mathcal{H}$.
\\ \\ {\bf Definition 1.} (see Theorem 2.10 in \cite{BB02}) Let $x\in\partial B_N$ be the Denjoy-Wolff point of $\varphi$ and $\lambda$ be the boundary dilatation coefficient of $\varphi$ at $x$. Define $$\mathbb{A}(\varphi):=\mbox{span}\{\nu\in\mathbb{C}^N: d\varphi_x(\nu)=\lambda\nu \ \mbox{and} <\nu, x>\ne 0\}$$
and $$\mathcal{A}\mathcal{G}(\varphi):=\bigcup\limits_{j=1}^{\infty}\mbox{ker} (d\varphi_x-\lambda I)^.$$
The spaces $\mathbb{A}(\varphi)$ and $\mathcal{A}\mathcal{G}(\varphi)$ are called the inner space and generalized inner space of $\varphi$ respectively.
\\ \par
Remark 2.11 in \cite{BB02} told us $\mathcal{A}\mathcal{G}(\varphi)$ is invariant for $d\varphi_x$ and $B_N\bigcap (\mathcal{A}\mathcal{G}(\varphi)+x)$ is the maximum (may be proper) invariant set of $\varphi$ in the ball. Any other invariant set of $\varphi$ is obtained as $B_N\bigcap (W+x)$ for $W\subset \mathcal{A}\mathcal{G}(\varphi)$ and $d\varphi_x(W)\subset W$.
\\ \\ {\bf Remark 2.} If $\mathbb{A}(\varphi)=\{0\}$, we see that $x$ is the only fixed point of $\varphi$. In this case, suppose that the restriction of $\varphi$ to $B_N\bigcap(W+x)$ is an automorphism for some  $W\subset \mathcal{A}\mathcal{G}(\varphi)$ and $d\varphi_x(W)\subset W$. Moreover,  if there exists an automorphism $\rho$ of $B_N$ that satisfies $\rho(B_N\bigcap(W+x))\subset B_k$ with $k=\mbox{dim}\, W>0$, such that the essential normality of $\varphi$ is not changed when conjugated by $\rho$. Then by Theorem 4.1,  $C_\varphi$ is not essentially normal on $\mathcal{H}$.
\\ \\ {\bf Definition 2.} (see  \cite{BCD07}) Let $z_0\in B_N$ be a fixed point of $\varphi$, we define the unitary space of $\varphi$ at $z_0$ by
$$L_U(\varphi, z_0):=\bigoplus\limits_{|\lambda|=1}\mbox{ker}(d\varphi_{z_0}-\lambda I)^N.$$
\\ {\bf Remark 3.} Write $L=L_U(\varphi, z_0)$. According to the proof of Theorem 3.2 (1) in \cite{BCD07}, if $k=\mbox{dim}\, L>0$, then
$\varphi(B_N\bigcap L)\subset B_N\bigcap L$ and  the restriction of $\varphi$ to $B_N\bigcap L$ is an automorphism. Now, suppose that $\varphi|_{B_N\bigcap L}$
 is not a rotation and the essential normality of $\varphi$ is not changed, when conjugated by an automorphism $\rho$ of $B_N$ with $\rho(B_N\bigcap L)\subset B_k$. Using Theorem 4.1, we obtain that $C_\varphi$ is not essentially normal on $\mathcal{H}$.
\\ \par
At last, we give another  class of linear fractional self-maps of $B_N$ whose  corresponding composition operators are not essentially normal on $\mathcal{H}$. The idea comes from the proofs of Proposition 2 in \cite{MW05} and Theorem 2.3 in \cite{JC14}.

\begin{theorem}\label{the 4.4} Let  $\varphi$ be  a linear fractional self-map   of $B_N$ with only one interior fixed point $z_0$ on $\overline{B_N}$.
 If $p=\mbox{dim}\,L_U(\varphi, z_0)=0$ and  $||\varphi||_{\infty}=1$, then $C_\varphi$ is not essentially normal on  $H^2(B_N)$ or $A^2_s(B_N)$  ($s>-1$).
\end{theorem}
\proof
 Let $\varphi^n$
denote the $n$-th iterate of $\varphi$. First, the proof of Theorem 2.3 in  \cite{JC14}  gives $||\varphi^n||_{\infty}<1$ for all $n\ge M$, where $M$ is a positive integer. Let $m=\min\{n:\ ||\varphi^n||_{\infty}<1 \}$. Since  $||\varphi||_{\infty}=1$, we see that $m>1$. This means that $||\varphi^m||_{\infty}<1$
and $||\varphi^{m-1}||_{\infty}=1$. Hence, there exist $\zeta, \eta\in\partial B_N$
such that $\varphi^{m-1}(\zeta)=\eta$. Applying Lemma 1 of \cite{MW05} repeatedly, we  see that $\sigma^{m-1}(\eta)=\zeta$. Which implies that $||\sigma^{m-1}||_{\infty}=1$. Moreover, it is easy to see $||\sigma^m||_{\infty}<1$, otherwise, we can use similar argument to obtain $||\varphi^m||_{\infty}=1$.

Now, we have $\sigma^{m-1}(\eta)=\zeta$ and $|\sigma^m(\eta)|<1$. Using Lemma 1 of \cite{MW05} again, we get $\varphi(\zeta)=\sigma^{m-2}(\eta)\in\partial B_N$.  For $p\in B_N$, let $k_p=K_p/||K_p||_{\mathcal{H}}$ be the normalized reproducing kernel on $\mathcal{H}$. Similar to the proof of Proposition 2 in \cite{MW05}, since $C_\varphi=(T_gC_\sigma T^*_h)^*=T_h C^*_\sigma T^*_g$, we have{\setlength\arraycolsep{2pt}
\begin{eqnarray*}  ||[C^*_\varphi, C_\varphi]k_p||_{\mathcal{H}}&\ge & |||C_\varphi k_p||_{\mathcal{H}}^2-||C^*_\varphi k_p||_{\mathcal{H}}^2|
\\ &=& \biggl|\frac{||\overline{g(p)}hK_{\sigma(p)}||_{\mathcal{H}}^2}{||K_p||_{\mathcal{H}}^2}-
\frac{||K_{\varphi(p)}||_{\mathcal{H}}^2}{||K_p||_{\mathcal{H}}^2}\biggr|
\\ &\ge & \biggl(\frac{1-|p|^2}{1-|\varphi(p)|^2}\biggr)^t-|g(p)|^2||h||^2_{\infty}\biggl(\frac{1-|p|^2}{1-|\sigma(p)|^2}\biggr)^t,
\end{eqnarray*}}where $t=N$ when $\mathcal{H}=H^2(B_N)$ and $t=N+s+1$ when  $\mathcal{H}=A_s^2(B_N)$ ($s>-1$). Choose a sequence $\{p_n\}$
 in $B_N$ with $p_n$ tending to $\zeta$ as $n\to\infty$, since $\varphi(\zeta)\in\partial B_N$ and $\sigma(\zeta)=\sigma^m(\eta)\in B_N$,  when  $n\to\infty$, we have  $$\biggl(\frac{1-|p_n|^2}{1-|\varphi(p_n)|^2}\biggr)^t\to d_\varphi(\zeta)^{-t}$$
 from Julia-Carath\'{e}odory Theorem in $B_N$ (see Theorem 2.8 of \cite{CM1}) and $$\biggl(\frac{1-|p_n|^2}{1-|\sigma(p_n)|^2}\biggr)^t\to 0,$$ where
 $d_\varphi(\zeta)=\liminf_{z\to\zeta}(1-|\varphi(z)|)/(1-|z|)>0$.  Hence, $[C^*_\varphi, C_\varphi]$ is not compact on  $\mathcal{H}$ and we obtain the desired conclusion. \ \ $\Box$
\\ \\ \\
\centerline{ APPENDIX: ESTIMATION FOR $\int_{\partial B_N}|<\zeta, a>^k<\zeta,\eta>^m|^2d\sigma(\zeta)$}
\\ \par We estimate the integral  $$\int_{\partial B_N}|<\zeta, a>^k<\zeta,\eta>^m|^2d\sigma(\zeta)$$ when $N=2$ and $k=2$.

First, let $m=1$. Since
 {\setlength\arraycolsep{2pt}
\begin{eqnarray*} && <\zeta, a>^2<\zeta,\eta>=(\zeta_1\overline{a_1}+\zeta_2\overline{a_2})^2(\zeta_1\overline{\eta_1}+\zeta_2\overline{\eta_1})
\\ &=& \zeta_1^3\overline{a_1}^2\overline{\eta_1}+2\zeta_1^2\zeta_2\overline{a_1}\overline{a_2}\overline{\eta_1}
+\zeta_1\zeta_2^2\overline{a_2}^2\overline{\eta_1}
+\zeta_1^2\zeta_2\overline{a_1}^2\overline{\eta_2}
+2\zeta_1\zeta_2^2\overline{a_1}\overline{a_2}\overline{\eta_2}+\zeta_2^3\overline{a_2}^2\overline{\eta_2}.
\end{eqnarray*}}Thus, {\setlength\arraycolsep{2pt}
\begin{eqnarray*} && \int_{\partial B_2}|<\zeta, a>^2<\zeta,\eta>|^2d\sigma(\zeta)
\\ &=&\int_{\partial B_2} \biggl|\zeta_1^3\overline{a_1}^2\overline{\eta_1}+2\zeta_1^2\zeta_2\overline{a_1}\overline{a_2}\overline{\eta_1}
+\zeta_1\zeta_2^2\overline{a_2}^2\overline{\eta_1}
+\zeta_1^2\zeta_2\overline{a_1}^2\overline{\eta_2}
+2\zeta_1\zeta_2^2\overline{a_1}\overline{a_2}\overline{\eta_2}+\zeta_2^3\overline{a_2}^2\overline{\eta_2}\biggr|^2d\sigma(\zeta)
\\ &=&\frac{(2-1)!1!}{(2-1+1+2)!}\biggl(|a_1|^4|\eta_1|^2\times 3!+4|a_1|^2|a_2|^2|\eta_1|^2\times 2!+|a_2|^4|\eta_1|^2\times 2!
\\ && +|a_1|^4|\eta_2|^2\times 2!+ 4 |a_1|^2|a_2|^2|\eta_2|^2\times 2! +|a_2|^4|\eta_2|^2\times 3!
 +2 |a_1|^2 a_1\overline{a_2}\overline{\eta_1}\eta_2\times 2!
 \\ && +2|a_2|^2 a_1\overline{a_2}\overline{\eta_1}\eta_2\times 2! +2|a_1|^2 \overline{a_1} a_2\eta_1\overline{\eta_2}\times 2!+2|a_2|^2 \overline{a_1} a_2\eta_1\overline{\eta_2}\times 2!
\biggr)
\\ &=&\frac{(2-1)!1!}{(2-1+1+2)!}\biggl[\biggl(4 |a_1|^4|\eta_1|^2+2 |a_1|^4|\eta_1|^2+4|a_1|^2|a_2|^2|\eta_1|^2+4|a_1|^2|a_2|^2|\eta_1|^2
\\ && +2|a_2|^4|\eta_1|^2+2|a_1|^4|\eta_2|^2+4 |a_1|^2|a_2|^2|\eta_2|^2+4 |a_1|^2|a_2|^2|\eta_2|^2+4|a_2|^4|\eta_2|^2+2|a_2|^4|\eta_2|^2
\\ &&  +4 |a_1|^2 a_1\overline{a_2}\overline{\eta_1}\eta_2
+4|a_2|^2 a_1\overline{a_2}\overline{\eta_1}\eta_2
 +4|a_1|^2 \overline{a_1} a_2\eta_1\overline{\eta_2}
+4|a_2|^2 \overline{a_1} a_2\eta_1\overline{\eta_2}
\\ &=&\frac{(2-1)!1!}{(2-1+1+2)!}\biggl[\biggl(4 |a_1|^4|\eta_1|^2+4|a_1|^2|a_2|^2|\eta_1|^2+4 |a_1|^2|a_2|^2|\eta_2|^2+4|a_2|^4|\eta_2|^2
\\ &&  +4 |a_1|^2 a_1\overline{a_2}\overline{\eta_1}\eta_2
+4|a_2|^2 a_1\overline{a_2}\overline{\eta_1}\eta_2
 +4|a_1|^2 \overline{a_1} a_2\eta_1\overline{\eta_2}
+4|a_2|^2 \overline{a_1} a_2\eta_1\overline{\eta_2}\biggr)
\\ && +2\biggl( |a_1|^4|\eta_1|^2+2|a_1|^2|a_2|^2|\eta_1|^2
 +|a_2|^4|\eta_1|^2
 \\ && +|a_1|^4|\eta_2|^2+2 |a_1|^2|a_2|^2|\eta_2|^2+|a_2|^4|\eta_2|^2
\biggr)\biggr]
\\ &=&\frac{(2-1)!1!}{(2-1+1+2)!}\biggl[\biggl(4 |a_1|^4|\eta_1|^2+4|a_1|^2|a_2|^2|\eta_1|^2+4 |a_1|^2|a_2|^2|\eta_2|^2+4|a_2|^4|\eta_2|^2
\\ &&  +4 |a_1|^2 a_1\overline{a_2}\overline{\eta_1}\eta_2
+4|a_2|^2 a_1\overline{a_2}\overline{\eta_1}\eta_2
 +4|a_1|^2 \overline{a_1} a_2\eta_1\overline{\eta_2}
+4|a_2|^2 \overline{a_1} a_2\eta_1\overline{\eta_2}\biggr)
\\ &&+2(|a_1|^2+|a_2|^2)^2(|\eta_1|^2+|\eta_2|^2)\biggr]
\end{eqnarray*}}Using the  inequality{\setlength\arraycolsep{2pt}
\begin{eqnarray*} && |a_1|^2|\eta_2|^2+ |a_2|^2|\eta_1|^2
\ge 2 |a_1||a_2||\eta_1||\eta_2|
\\ && \ge 2\mbox{Re}\,(\overline{a_1} a_2\eta_1 \overline{\eta_2})
 =\overline{a_1} a_2\eta_1 \overline{\eta_2}+ a_1 \overline{a_2}\, \overline{\eta_1}\eta_2,
\end{eqnarray*}}We obtain {\setlength\arraycolsep{2pt}
\begin{eqnarray*} && \int_{\partial B_2}|<\zeta, a>^2<\zeta,\eta>|^2d\sigma(\zeta)
\\ &=&\frac{(2-1)!1!}{(2-1+1+2)!}\biggl[\biggl(4 |a_1|^4|\eta_1|^2+4|a_1|^2|a_2|^2|\eta_1|^2+4 |a_1|^2|a_2|^2|\eta_2|^2+4|a_2|^4|\eta_2|^2
\\ &&  +4 |a_1|^2( a_1\overline{a_2}\overline{\eta_1}\eta_2+ \overline{a_1} a_2\eta_1\overline{\eta_2})
+4|a_2|^2 (a_1\overline{a_2}\overline{\eta_1}\eta_2
 +\overline{a_1} a_2\eta_1\overline{\eta_2})\biggr)
+2|a|^4|\eta|^2\biggr]
\\ &\le &\frac{(2-1)!1!}{(2-1+1+2)!}\biggl[\biggl(4 |a_1|^4|\eta_1|^2+4|a_1|^2|a_2|^2|\eta_1|^2+4 |a_1|^2|a_2|^2|\eta_2|^2+4|a_2|^4|\eta_2|^2
\\ &&  +4 |a_1|^2(|a_1|^2|\eta_2|^2+ |a_2|^2|\eta_1|^2)
+4|a_2|^2 (|a_1|^2|\eta_2|^2+ |a_2|^2|\eta_1|^2)\biggr)
+2|a|^4|\eta|^2\biggr]
\\ &= &\frac{(2-1)!1!}{(2-1+1+2)!}\biggl[4 \biggl(|a_1|^4|\eta_1|^2+2|a_1|^2|a_2|^2|\eta_1|^2+ |a_2|^4|\eta_1|^2
\\ && + |a_1|^4|\eta_2|^2+ 2|a_1|^2|a_2|^2|\eta_2|^2+|a_2|^4|\eta_2|^2
  \biggr)
+2|a|^4|\eta|^2\biggr]
\\ &= &\frac{(2-1)!1!}{(2-1+1+2)!}(4|a|^4|\eta|^2+2|a|^4|\eta|^2)
\\&=& \frac{(2-1)!1!}{(2-1+1+2)!}(1+1)(1+2)|a|^4|\eta|^2
\\ &= &\frac{(2-1)!(1+2)!}{(2-1+1+2)!}|a|^4|\eta|^2=\frac{(2-1)!(1+2)!}{(2-1+1+2)!}|a|^{2\times 2}.
\end{eqnarray*}}

For $m\ge 2$, we have{\setlength\arraycolsep{2pt}
\begin{eqnarray*} && <\zeta, a>^2<\zeta,\eta>^m=(\zeta_1\overline{a_1}+\zeta_2\overline{a_2})^2(\zeta_1\overline{\eta_1}+\zeta_2\overline{\eta_1})^m
\\ &=& (\overline{a_1}^2\zeta_1^2+2\overline{a_1}\overline{a_2}\zeta_1\zeta_2+\overline{a_2}^2\zeta_2^2)
\biggl(\sum\limits_{|\alpha|=m}\frac{m!}{\alpha!}\overline{\eta_1}^{\alpha_1}\overline{\eta_2}^{\alpha_2}  \zeta_1^{\alpha_1} \zeta_2^{\alpha_2}\biggr)
\\ &=&
\sum\limits_{|\alpha|=m}\frac{m!}{\alpha!}\overline{a_1}^2\overline{\eta_1}^{\alpha_1}\overline{\eta_2}^{\alpha_2}  \zeta_1^{\alpha_1+2} \zeta_2^{\alpha_2}
+2\sum\limits_{|\alpha|=m}\frac{m!}{\alpha!}\overline{a_1}\overline{a_2}\overline{\eta_1}^{\alpha_1}\overline{\eta_2}^{\alpha_2}  \zeta_1^{\alpha_1+1} \zeta_2^{\alpha_2+1}
\\ && +\sum\limits_{|\alpha|=m}\frac{m!}{\alpha!}\overline{a_2}^2\overline{\eta_1}^{\alpha_1}\overline{\eta_2}^{\alpha_2}  \zeta_1^{\alpha_1} \zeta_2^{\alpha_2+2}.
\end{eqnarray*}}Hence,  {\setlength\arraycolsep{2pt}
\begin{eqnarray*} && \int_{\partial B_2}|<\zeta, a>^2<\zeta,\eta>^m|^2d\sigma(\zeta)
\\ &=&\int_{\partial B_2}\biggl|\sum\limits_{|\alpha|=m}\frac{m!}{\alpha!}\overline{a_1}^2\overline{\eta_1}^{\alpha_1}\overline{\eta_2}^{\alpha_2}  \zeta_1^{\alpha_1+2} \zeta_2^{\alpha_2}
+2\sum\limits_{|\alpha|=m}\frac{m!}{\alpha!}\overline{a_1}\overline{a_2}\overline{\eta_1}^{\alpha_1}\overline{\eta_2}^{\alpha_2}  \zeta_1^{\alpha_1+1} \zeta_2^{\alpha_2+1}
\\ && +\sum\limits_{|\alpha|=m}\frac{m!}{\alpha!}\overline{a_2}^2\overline{\eta_1}^{\alpha_1}\overline{\eta_2}^{\alpha_2}  \zeta_1^{\alpha_1} \zeta_2^{\alpha_2+2}\biggr|^2d\sigma(\zeta)
\\ &=& \frac{(2-1)!m!}{(2-1+m+2)!}\biggl(
\sum\limits_{|\alpha|=m}\frac{m!}{\alpha!}\frac{(\alpha_1+2)!\alpha_2!}{\alpha!}|a_1|^4|\eta_1|^{2\alpha_1}|\eta_2|^{2\alpha_2}
\end{eqnarray*}}
 {\setlength\arraycolsep{2pt}
\begin{eqnarray*}&& +4\sum\limits_{|\alpha|=m}\frac{m!}{\alpha!}\frac{(\alpha_1+1)!(\alpha_2+1)!}{\alpha!}|a_1|^2|a_2|^2|\eta_1|^{2\alpha_1}
|\eta_2|^{2\alpha_2}
+\sum\limits_{|\alpha|=m}\frac{m!}{\alpha!}\frac{\alpha_1!(\alpha_2+2)!}{\alpha!}
\\&& \times|a_2|^4|\eta_1|^{2\alpha_1}|\eta_2|^{2\alpha_2}
 +\sum\limits_{\gamma\ne\delta}\sum\limits_{\gamma+\alpha\ge \delta}\frac{2!}{\gamma!}\frac{2!}{\delta!}\frac{m!}{\alpha!}\frac{(\gamma+\alpha)!}{(\gamma+\alpha-\delta)!}
\overline{a}^\gamma a^\delta\overline{\eta}^\alpha \eta^{\gamma+\alpha-\delta} \biggr).
\end{eqnarray*}}We use the decomposition   in Lemma 3.3 to compute $$\frac{(\alpha_1+2)!}{\alpha_1!}\qquad \mbox{and}\qquad \frac{(\alpha_2+2)!}{\alpha_2!}.$$
Note that $$\frac{(\gamma+\alpha)!}{(\gamma+\alpha-\delta)!}$$  can also be decomposed  by using Lemma 3.3, but it is redundant to display, so we omit this part. Hence, {\setlength\arraycolsep{2pt}
\begin{eqnarray*} && \int_{\partial B_2}|<\zeta, a>^2<\zeta,\eta>^m|^2d\sigma(\zeta)
\\ &=&\frac{(2-1)!m!}{(2-1+m+2)!}\biggl(
\sum\limits_{|\alpha|=m}\frac{m!}{\alpha!}[\alpha_1(\alpha_1-1)+4\alpha_1+2]|a_1|^4|\eta_1|^{2\alpha_1}|\eta_2|^{2\alpha_2}
\\&& +4\sum\limits_{|\alpha|=m}\frac{m!}{\alpha!}(\alpha_1\alpha_2+\alpha_1+\alpha_2+1)|a_1|^2|a_2|^2|\eta_1|^{2\alpha_1}
|\eta_2|^{2\alpha_2}
\\ &&  +\sum\limits_{|\alpha|=m}\frac{m!}{\alpha!}[\alpha_2(\alpha_2-1)+4\alpha_2+2]|a_2|^4|\eta_1|^{2\alpha_1}|\eta_2|^{2\alpha_2}
\\ && + \sum\limits_{\gamma\ne\delta}\sum\limits_{\gamma+\alpha\ge \delta}\frac{2!}{\gamma!}\frac{2!}{\delta!}\frac{m!}{\alpha!}\frac{(\gamma+\alpha)!}{(\gamma+\alpha-\delta)!}
\overline{a}^\gamma a^\delta\overline{\eta}^\alpha \eta^{\gamma+\alpha-\delta} \biggr)
\\ &=&\frac{(2-1)!m!}{(2-1+m+2)!}\biggl(
m(m-1)|a_1|^4|\eta_1|^4\sum\limits_{\substack{|\alpha|=m \\ \alpha_1\ge 2}}\frac{(m-2)!}{(\alpha_1-2)!\alpha_2!}|\eta_1|^{2(\alpha_1-2)}|\eta_2|^{2\alpha_2}
\\&& +4m|a_1|^4|\eta_1|^2\sum\limits_{\substack{|\alpha|=m \\ \alpha_1\ge 1}}\frac{(m-1)!}{(\alpha_1-1)!\alpha_2!}|\eta_1|^{2(\alpha_1-1)}|\eta_2|^{2\alpha_2}
+2|a_1|^4\sum\limits_{|\alpha|=m}\frac{m!}{\alpha!}|\eta_1|^{2\alpha_1}|\eta_2|^{2\alpha_2}
\\&& +4m(m-1)|a_1|^2|a_2|^2|\eta_1|^2|\eta_2|^2\sum\limits_{\substack{|\alpha|=m \\ \alpha_1\ge 1,\  \alpha_2\ge 1}}\frac{(m-2)!}{(\alpha_1-1)!(\alpha_2-1)!}|\eta_1|^{2(\alpha_1-1)}|\eta_2|^{2(\alpha_2-1)}
\\&& +4m|a_1|^2|a_2|^2|\eta_1|^2\sum\limits_{\substack{|\alpha|=m \\ \alpha_1\ge 1}}\frac{(m-1)!}{(\alpha_1-1)!\alpha_2!}|\eta_1|^{2(\alpha_1-1)}|\eta_2|^{2\alpha_2}
\\&& +4m|a_1|^2|a_2|^2|\eta_2|^2\sum\limits_{\substack{|\alpha|=m \\ \alpha_2\ge 1}}\frac{(m-1)!}{\alpha_1!(\alpha_2-1)!}|\eta_1|^{2\alpha_1}|\eta_2|^{2(\alpha_2-1)}
\\&& +4|a_1|^2|a_2|^2\sum\limits_{|\alpha|=m}\frac{m!}{\alpha!}|\eta_1|^{2\alpha_1}
|\eta_2|^{2\alpha_2}
  +m(m-1)|a_2|^4|\eta_2|^4\sum\limits_{\substack{|\alpha|=m \\ \alpha_2\ge 2}}\frac{(m-2)!}{\alpha_1!(\alpha_2-2)!}|\eta_1|^{2\alpha_1}|\eta_2|^{2(\alpha_2-2)}
 \end{eqnarray*}}
{\setlength\arraycolsep{2pt}
\begin{eqnarray*} && +4m|a_2|^4|\eta_2|^2\sum\limits_{\substack{|\alpha|=m \\ \alpha_2\ge 1}}\frac{(m-1)!}{\alpha_1!(\alpha_2-1)!}|\eta_1|^{2\alpha_1}|\eta_2|^{2(\alpha_1-1)}
+2|a_2|^4\sum\limits_{|\alpha|=m}\frac{m!}{\alpha!}|\eta_1|^{2\alpha_1}|\eta_2|^{2\alpha_2}
\\&& +\sum\limits_{\gamma\ne\delta}\sum\limits_{\gamma+\alpha\ge \delta}\frac{2!}{\gamma!}\frac{2!}{\delta!}\frac{m!}{\alpha!}\frac{(\gamma+\alpha)!}{(\gamma+\alpha-\delta)!}
\overline{a}^\gamma a^\delta\overline{\eta}^\alpha \eta^{\gamma+\alpha-\delta} \biggr).
\end{eqnarray*}}Similar to the idea when $m=1$, we use the inequality  $$|a_1|^2|\eta_2|^2+ |a_2|^2|\eta_1|^2
 \ge \overline{a_1} a_2\eta_1 \overline{\eta_2}+ a_1 \overline{a_2}\, \overline{\eta_1}\eta_2$$
 and $$|a_1|^4|\eta_2|^4+ |a_2|^4|\eta_1|^4
 \ge \overline{a_1}^2 a_2^2\eta_1^2 \overline{\eta_2}^2+ a_1^2 \overline{a_2}^2\, \overline{\eta_1}^2\eta_2^2$$  to estimate the last sum above.
Therefore, {\setlength\arraycolsep{2pt}
\begin{eqnarray*} && \int_{\partial B_2}|<\zeta, a>^2<\zeta,\eta>^m|^2d\sigma(\zeta)
\\ &\le & \frac{(2-1)!m!}{(2-1+m+2)!}\biggl[m(m-1)\biggl(|a_1|^4|\eta_1|^4+2|a_1|^4|\eta_1|^2|\eta_2|^2+|a_1|^4|\eta_2|^4
\\&& +2|a_1|^2|a_2|^2|\eta_1|^4+ 4|a_1|^2|a_2|^2|\eta_1|^2|\eta_2|^2+2|a_1|^2|a_2|^2|\eta_2|^4
\\ && +|a_2|^4|\eta_1|^4+2|a_2|^4|\eta_1|^2|\eta_2|^2+|a_2|^4|\eta_2|^4\biggr)
 \times(|\eta_1|^2+|\eta_2|^2)^{m-2}
\\ && +4m\biggl(|a_1|^4|\eta_1|^2+|a_1|^4|\eta_2|^2+2|a_1|^2|a_2|^2|\eta_1|^2+2|a_1|^2|a_2|^2|\eta_2|^2+|a_2|^4|\eta_1|^2+|a_2|^4|\eta_2|^2\biggr)
\\ && \times (|\eta_1|^2+|\eta_2|^2)^{m-1}
 + 2(|a_1|^4+2|a_1|^2|a_2|^2+|a_2|^4)(|\eta_1|^2+|\eta_2|^2)^m\biggr]
\\ &=& \frac{(2-1)!m!}{(2-1+m+2)!}\biggl[m(m-1)(|a_1|^2+|a_2|^2)^2(|\eta_1|^2+|\eta_2|^2)^2(|\eta_1|^2+|\eta_2|^2)^{m-2}
\\ && +4m(|a_1|^2+|a_2|^2)^2(|\eta_1|^2+|\eta_2|^2)(|\eta_1|^2+|\eta_2|^2)^{m-1}+2(|a_1|^2+|a_2|^2)^2(|\eta_1|^2+|\eta_2|^2)^m\biggr]
\\ &=& \frac{(2-1)!m!}{(2-1+m+2)!}\biggl[m(m-1)|a|^4|\eta|^{2m}+4m|a|^4|\eta|^{2m}+2|a|^4|\eta|^{2m}\biggr]
\\ &=&  \frac{(2-1)!m!}{(2-1+m+2)!}(m+1)(m+2)|a|^4|\eta|^{2m}
\\ &=&  \frac{(2-1)!(m+2)!}{(2-1+m+2)!}|a|^{2\times 2}.
\end{eqnarray*}}For more general $N$ and $k$, we can use similar idea to compute the previous integral.
\\ \\ \\
\centerline{ ACKNOWLEDGEMENTS }
\\ \par This paper was written during the first author's visiting at The College at Brockport, State University of New York. SUNY Brockport provided an excellent environment for working on this paper. Shanghai Municipal Education Commission  provided the financial support during her visiting. She would like to express her gratitude to them. She is also grateful to  Matthew A. Pons for providing the references \cite{MP06} and \cite{Po11}.

\bibliographystyle{amsplain}

\end{document}